\documentclass[preprint,12pt]{elsarticle}

\usepackage{graphicx}

\usepackage{hyperref}

\usepackage{multirow}

\usepackage[ruled,lined,algonl,boxed]{algorithm2e}

\usepackage{amsmath,amssymb}

\usepackage{amsthm}

\usepackage[utf8]{inputenc}
\usepackage{float}
\usepackage[table,xcdraw,dvipsnames]{xcolor}
\usepackage{hhline}
\usepackage{makecell}
\graphicspath{ {./Graphs/} }
\usepackage{pgfplots}
\pgfplotsset{width=150mm,height=75mm,scaled y ticks=false,compat=1.9}

\allowdisplaybreaks

\newtheorem{thm}{Theorem}
\newtheorem{cor}[thm]{Corollary}
\newdefinition{rmk}{Remark}
\newproof{pf}{Proof}
\newproof{pot}{Proof of Theorem \ref{thm2}}

\journal{}

\usepackage{amsopn}

\begin{document}

\begin{frontmatter}



\title{Weighted domination models and randomized heuristics}


\author[label2]{Lukas Dijkstra}
\ead{dijkstral@cardiff.ac.uk}
\author[label2]{Andrei Gagarin}
\ead{gagarina@cardiff.ac.uk}
\author[label2a]{Vadim Zverovich}
\ead{vadim.zverovich@uwe.ac.uk}
\address[label2]{School of Mathematics, Cardiff University, Cardiff, UK}
\address[label2a]{\,\href{https://www.uwe.ac.uk/research/centres-and-groups/mathematics-and-statistics}{Mathematics and Statistics Research Group}, University of the West of England,\\ Bristol, UK}


\begin{abstract}
We consider the minimum weight and smallest weight minimum-size dominating set problems in vertex-weighted graphs and networks. The latter problem is a two-objective optimization problem, which is different from the classic minimum weight dominating set problem that requires finding a dominating set of the smallest weight in a graph without trying to optimize its cardinality. 
In other words, the objective of minimizing the size of the dominating set in the two-objective problem can be considered as a constraint, i.e. a particular case of finding Pareto-optimal solutions.
First, we show how to reduce the two-objective optimization problem to the minimum weight dominating set problem by using Integer Linear Programming formulations.
Then, under different assumptions, the probabilistic method is applied to obtain upper bounds on the minimum weight dominating sets in graphs. 
The corresponding randomized algorithms for finding small-weight dominating sets in graphs are described as well. 
Computational experiments are used to illustrate the results for two different types of random graphs.
\end{abstract}


\end{frontmatter}


\section{Introduction}
\label{Intro}
\subsection{Basic notions, notation, and motivation}
\label{notation}

We consider undirected simple vertex-weighted graphs 
$$G=(V,E,w\!: V\rightarrow \mathbb{R}),$$ 
where $V=V(G)=\{v_1,v_2,...,v_n\}$ is the set of vertices, $E=E(G)=\{e_1,e_2,...,e_m\}$ is the set of edges of $G$, 
and $w\!: V\rightarrow \mathbb{R}$ is a weight\slash cost function 
that assigns a certain weight $w_i=w(v_i)$
to each vertex $v_i$ of $G$, $i=1,...,n$. 
The neighbourhood of a vertex $v$ in $G$ is denoted by $N(v)$, i.e. $N(v)=\{u\ |\ vu\in E,\ u\not= v\}$. Any vertex in $N(v)$ is a neighbour of $v$.
The closed neighbourhood of $v$ is denoted by $N[v] = N(v)\cup \{v\}$. 
For a set of vertices $A\subseteq V$, 
$$
N(A)=\bigcup_{v\in A} N(v) \quad \mbox{and} \quad N[A]= N(A)\cup A.
$$
The degree of a vertex $v_i$ is the number of its neighbours  and is denoted by $d_i=d(v_i)$, $i=1,...,n$. 
A sequence of vertex degrees of $G$ is denoted by $\bar{d}=(d_1,d_2,...,d_n)$. 
The minimum and maximum vertex degrees of $G$ are denoted by $\delta=\delta(G)$ and $\Delta=\Delta(G)$, respectively. 

A subset $X\subseteq V(G)$ is
called a {\it dominating set} of $G$ if every vertex not in $X$ is
adjacent to at least one vertex in $X$. 
The minimum cardinality of
a dominating set of $G$ is called the {\it domination number} of $G$ and denoted by 
$\gamma(G)$.
We denote by $\gamma_w(G)$ the smallest weight of a dominating set in a vertex-weighted graph $G$, 
and by $\gamma_w^*(G)$ the smallest weight of a minimum-cardinality dominating set $D$ in $G$. 
Clearly,  $\displaystyle \gamma_w(G)\le \gamma_w^*(G)$. 

The total weight of the graph is
$$
w_G=\sum_{v_i\in V} w_i.
\vspace*{-0.3cm}$$
Also,
$$
\displaystyle w_{\rm max}=\max_{1\le i \le n} w_i,\quad 
\displaystyle w_{\rm min}=\min_{1\le i \le n} w_i, \quad  \displaystyle w_{\rm ave}=\frac{w_G}{n},
$$
so that
 $$\displaystyle w_{\rm max}\ge w_{\rm ave}\ge w_{\rm min}.$$ 

Weighted domination in graphs and networks can be used, for example, for modelling a problem  of the placement of a small number of transmitters in a communication network such that every site in the network either has a transmitter or is connected by a direct communication link to a site that has such a transmitter. 
In addition, there are  some `costs' associated with placing a transmitter in each particular location of the network (i.e. a vertex of the corresponding graph). 
The minimum weight  dominating set problem usually does not place any restrictions on the size of the dominating set, i.e. the number of transmitters in this case -- it only requires us to find a smallest weight\slash cost dominating set in a vertex-weighted graph. 
However, the total emitted radiation in the environment would be smaller
with fewer transmitters installed.
Therefore, if the weight associated with each vertex of the graph is considered as a cost function or some nuisance measurement parameter (e.g. a level of noise at the site of an installed transmitter), the corresponding problem becomes a multi-criteria minimization problem, where the objective is to find 
a minimum cost smallest-cardinality dominating set in a weighted graph.

The weights of vertices in the graph can also indicate a local impact or a level of influence of placing a facility in each particular location; for example, see positive influence dominating sets in social networks \cite{GPC2014,TCS2011}. 
In this case, we have a multi-criteria optimization problem of finding 
a smallest-cardinality dominating set $X$ in a graph $G$
because of the limited availability of the resource. 
On the other hand, the set $X$ needs to provide the maximum positive level of impact on the whole network represented by $G$. 
In other words, it is necessary to find a smallest-cardinality dominating set $X$ in $G$ such that its total weight $w(X)=\sum_{v_i\in X} w_i$ is maximized.
More examples of applications and studies of dominating sets for modelling real-life problems in networks can be found in \cite{Zve2, GC2018, CH1977}.

Assuming that all weights $w_i$ are non-negative and $w_{\rm max} > 0$, the last problem can be reduced to the {corresponding} minimization problem as follows. 
We can scale the weights at vertices of $G$ to the interval $[0,1]$, 
for example by dividing each of them by a positive constant $w_{\rm max}$, 
so that $\alpha w_i\in [0,1]$, where $\alpha=1/w_{\rm max}$.
Then, replace all original weights $w_i$ by the values $\psi_i=1-\alpha w_i$.
Clearly, $\psi_i\in [0,1]$ for $i=1,...,n$. 
Now, the problem of finding a minimum-size dominating set $X$ of the largest possible weight $w(X)=\sum_{v_i\in X} w_i$ in $G$ is equivalent to maximizing $\alpha \sum_{v_i\in X} w_i$ on the set of all minimum-size 
dominating sets $X\subseteq V(G)$.
Assuming that  
$\gamma(G)$ is known,
this is equivalent to minimizing the sum $\sum_{v_i\in X} \psi_i$ with an additional constraint on the cardinality of the dominating set, i.e. $|X|=\gamma(G)$, and we minimize
$$
\sum_{v_i\in X} \psi_i=|X| - \alpha \sum_{v_i\in X} w_i= \gamma(G) - \alpha \sum_{v_i\in X} w_i,
$$ 
where $\gamma(G)$ and $\alpha$ are constants.
Notice that $$w(X)=w_{\rm max}\cdot\left(\gamma(G) - \sum_{v_i\in X} \psi_i\right).$$ 
Therefore, we can focus on the minimization version of the general two-objective optimization problem.
{Notice that, in the classic single-objective minimum weight dominating set problem, vertices of weight zero and their neighbours can be removed from the graph.}

\subsection{Related results and complexity issues}
\label{related}

The problem of finding the exact value of $\gamma(G)$ in a graph $G$ is one of the classic NP-hard problems \cite{GJ1979}. Moreover, the problem is known to be APX-hard (e.g. see \cite{RS1997}) and not fixed parameter tractable \cite{DF1999}. 
Hence, the problem of finding a minimum or 
maximum weight smallest-cardinality dominating set in a graph $G$ is also NP-hard and not fixed parameter tractable in general. 
Therefore, it is necessary to have efficient and effective heuristic algorithms and methods for finding small-size and light- or heavy-weight dominating sets in graphs. 
Also, to estimate the quality of a given dominating set,
it is important to have good bounds for the domination number $\gamma(G)$ and for the smallest or largest weight $w(D)$ of  dominating sets $D\subseteq V(G)$.
The following upper bounds for the domination number $\gamma(G)$, which can be obtained using the probabilistic method, are well-known.


\begin{thm} [\cite{AS1992,HHS1998}] \label{arn}
For any graph $G$,
\begin{equation}\label{classical}
\gamma(G) \le {\ln(\delta+1)+1 \over \delta+1} n.
\end{equation}
\end{thm}


\begin{thm}[\cite{CR1997,HHS1998,HPV1999,Lee1998}] \label{Caro} 
For any graph $G$ with $\delta\ge 1$,
\begin{equation}\label{Car-Rod}
\gamma(G) \le \left(1-{\delta \over
{(1+\delta)^{1+1/\delta}}}\right) n.
\end{equation}
\end{thm}


These upper bounds are known to be asymptotically sharp by using Alon's construction (e.g. see \cite{ZP2011}).
Also, by applying the probabilistic method, these upper bounds have been generalized for directed graphs \cite{Lee1998} and for several other domination parameters in simple unweighted graphs  \cite{GPZ2009,GPZ2013,Zve2}.

Chen et al. \cite{CI2004} showed that the single-objective minimum weight dominating set problem is APX-hard and provided a simple greedy heuristic achieving an $O(\log n)$ approximation ratio, which is asymptotically the best possible in that case. 
The authors of \cite{CI2004} also gave a randomized rounding heuristic for a linear relaxation of an Integer Linear Programming (ILP) formulation of the problem. 
The authors in \cite{GPC2014} used a randomized rounding heuristic for a linear relaxation of an ILP formulation for a variation of this problem.
Notice that the minimum size dominating set problem can be considered as a particular case of the minimum weight dominating set problem by assigning unit weights to all vertices of the graph.
It is easy to see that the (two-objective) smallest or largest weight minimum-size dominating set problem in vertex-weighted graphs is APX-hard as a generalization of the minimum-size dominating set problem in simple graphs.

Also, notice that the domination number $\gamma(G)$ of a graph $G$ exhibits the asymptotic property of having two points of concentration in random graphs \cite{Weber1981} (and digraphs \cite{Lee1998}).
However, the asymptotic property of a ``typical" graph stated in this kind of theorems cannot be used as a bound for the domination number $\gamma(G)$ of a given graph $G$, 
it does not help to determine the domination number exactly, 
and it does not provide any ideas how to find the corresponding dominating sets of size $\gamma(G)$ in $G$.

In this paper, we focus on the (two-objective) smallest or largest weight minimum-size dominating set problem and 
direct applications of the probabilistic method to tackle the related (single-objective) minimum weight dominating set problem in vertex-weighted graphs. This allows us to analyze deterministic and  several heuristic solution methods 
for these two problems 
by considering corresponding computational results.
First, in Section \ref{sec:deterministic-heuristic}, we consider an ILP formulation of the two-objective optimization problem and show its connection and reduction to the single-objective minimum weight dominating set problem.
Then, probabilistic constructions, the corresponding randomized algorithms, and upper bounds are described for the single-objective minimum weight dominating set problem in Section \ref{Section2}. 
The new upper bounds presented in Section \ref{Section2} can be considered as a generalization of the classic upper bounds (\ref{Car-Rod}) and (\ref{classical}).
To illustrate the concepts and better analyze the results, some computational experiments with random graphs are described in Section \ref{sec:experiments}. 
Section \ref{sec:conclusion} provides a summary of our findings and conclusions.


\section{ILP formulation and reduction to the minimum weight problem}
\label{sec:deterministic-heuristic}

For simplicity, we focus on the minimization version of the problems with positive weights $w_i>0$ for all vertices of $G$ by trying to find a {minimum-size} dominating set $X$ in a vertex-weighted graph $G$ such that $w(X)$ is the smallest possible.
Generic ILP solvers are frequently used to solve applied optimization problems in industry. 
In this section, we describe how to formulate the two-objective minimum weight smallest-cardinality dominating set problem as an ILP problem to solve it deterministically. We also show connections with the single-objective minimum weight dominating set problem and provide a reduction to the latter problem.

Given a graph $G$, the problem of finding the exact value of $\gamma_w^*(G)$ and the corresponding dominating set $X\subseteq V(G)$ can be formulated as an ILP problem as follows. A $(0,1)$-decision variable $x_i\in \{0,1\}$ is associated with each vertex $v_i\in G$ to indicate whether the vertex is in the solution set $X$ or not, i.e. $x_i=1$ if and only if $v_i$ is in the smallest weight min-size dominating set $X$ of $G$, otherwise $x_i=0$, $i=1,...,n$. 
Here is the ILP problem formulation:\\
\begin{equation}
\begin{array}{ll@{}lrl}

\text{minimize}  & z(x_1,x_2,\ldots,x_n)=\displaystyle\sum\limits_{i=1}^{n} x_i + \sum\limits_{i=1}^{n} \frac{w_i}{w_G}&x_{i} &\\[7mm]

\text{subject to:}& \displaystyle\sum\limits_{v_{i} \in N[v_j]} x_{i} \geq 1,  &j=1 ,\dots, n, & \\[7mm]

                 &                                                x_{i} \in \{0,1\}, &i=1 ,\dots,n. &
\end{array}
\label{ILP1}
\end{equation}

Assuming $G$ is non-trivial and non-empty, we have $0< \frac{w_i}{w_G}<1$, $i=1,...,n$, and  $0< \sum\limits_{i=1}^{n} \frac{w_i}{w_G} x_{i}< 1$ in the objective function {for any non-trivial feasible solution} of problem $(\ref{ILP1})$. 
Therefore, at optimum, $z^*=z(x_1^*,x_2^*,...,x_n^*)$, we have $\displaystyle \gamma(G)=\sum\limits_{i=1}^{n} x_i^*$ and $\displaystyle \gamma_w^*(G)=w_G \cdot \sum\limits_{i=1}^{n} \frac{w_i}{w_G} x_{i}^* = \sum\limits_{i=1}^{n} w_i x_i^*$, i.e. $\displaystyle \gamma_w^*(G) = w_G\cdot (z^* - \gamma(G))$. Reassigning the graph vertex weights to $\displaystyle w'_i=1+\frac{w_i}{w_G}$, $i=1,...,n$, the ILP formulation $(\ref{ILP1})$ becomes 
\begin{equation}
\begin{array}{ll@{}lrl}

\text{minimize}  & z(x_1,x_2,\ldots,x_n)=\displaystyle\sum\limits_{i=1}^{n} \left(1+ \frac{w_i}{w_G}\right)&x_{i} =\displaystyle\sum\limits_{i=1}^{n} w'_i x_{i} & \\[7mm]

\text{subject to:}& \displaystyle\sum\limits_{v_{i} \in N[v_j]} x_{i} \geq 1,  &j=1 ,\dots, n, & \\[7mm]

                 &                                                x_{i} \in \{0,1\}, &i=1 ,\dots,n. &
\end{array}
\label{ILP2}
\end{equation}
This is the single-objective optimization problem of finding $\gamma_{w'}(G)$
and the corresponding minimum weight  dominating set in $G$ with respect to the vertex weights $w'_i$, $i=1,...,n$. 
Clearly, $\gamma_{w'}(G)=z^*=z(x_1^*,x_2^*,...,x_n^*)$ at optimum in (\ref{ILP1}) and (\ref{ILP2}), and $\gamma_w^*(G) = w_G\cdot (\gamma_{w'}(G) - \gamma(G))$, as above.
In other words, this provides a reduction from the two-objective problem of finding $\gamma_w^*(G)$ in $G$ to the single-objective problem of finding $\gamma_{w'}(G)$ in $G$. 

Notice that here problem $(\ref{ILP2})$ is a special case of the general minimum weight dominating set problem because the new weights are restricted to be $1< w'_i <2$, $i=1,...,n$. This may allow us to better understand the initial problem $(\ref{ILP1})$ and help to find more efficient solution methods.
{On the other hand, for the general minimum weight dominating set problem, we can always assume that the vertex weights belong to the interval $(0,1)$; for example, we can remove zero-weight vertices and their neighbours from the graph and divide the remaining weights by $w_{max}+1$.}
Also, notice that for the linear relaxation of problems (\ref{ILP1}) and (\ref{ILP2}) with $x_{i} \in [0,1]$, $i=1,...,n$, the objective function value $z'=z(x'_1,x'_2,...,x'_n)$ at the linear relaxation optimum provides lower bounds on $\gamma(G)$ and $\gamma_{w'}(G)$, i.e. $\gamma(G)\ge \lfloor z'\rfloor$ and $\gamma_{w'}(G)\ge z'$. 
{For $\gamma^*_w(G)$, more complicated lower bounds can be deduced. For instance, if $\gamma(G)\le U$, where $U$ is an upper bound, we get $\gamma^*_w(G)\ge w_G (z' - U)$, where we can trivially take $U=\sum_{i=1}^{n} \lceil x'_i\rceil$.}
On the other hand, any feasible solution to ILP problems (\ref{ILP1}) and (\ref{ILP2}), e.g. randomized rounding of their linear relaxation, allows us to recover the dominating set size and its weight.
Therefore, to tackle the problem of finding $\gamma_w^*(G)$ formulated in (\ref{ILP1}), one can use general tools to solve the minimum weight dominating set problem in a graph.



\section{Probabilistic constructions, upper bounds, and\\ randomized heuristics}
\label{Section2}

In this section, we focus on the minimum weight dominating set problem, i.e. finding upper bounds for $\gamma_w(G)$ in a graph $G$ analytically and algorithmically. Notice that, in view of the problem reduction in Section \ref{sec:deterministic-heuristic}, the other two similar problems for $\gamma(G)$ and $\gamma_w^*(G)$ in $G$ can be considered as particular cases of the above more general problem.

In the case of optimization of dominating sets in vertex-weighted graphs, it is reasonable to consider some vertex degrees,
e.g. the minimum, mean, or median vertex degree of the graph \cite{Zve2,GC2018}, and vertex weights as key parameters that influence the likelihood of inclusion of a vertex into an optimal dominating set with respect to its size and weight. 
For example, suppose we deal with the maximization version of the problem searching for a reasonably small-size dominating set $X$ of the maximum total weight $w(X)$ in a graph $G$. 
Then, it is plausible that a vertex of high degree and heavy weight is more likely to be included into $X$ than a vertex of lower degree and  light weight.
We will use this and similar assumptions in our applications of the probabilistic method and the corresponding randomized algorithms and techniques. 
Although  in this case, there is clearly a certain trade-off between the weight and degree parameters for a vertex to be included into an optimal dominating set of $G$.

As above, we assume positive weights $w_i>0$ for all vertices  
in a given vertex-weighted graph $G$, $i=1,2,...,n$, and
try  to find a dominating set $X$ in $G$ such that $w(X)$ is the smallest possible.
First, we generalize the probabilistic approach that was used to obtain the classic upper bounds  of Theorems \ref{arn} and \ref{Caro} for the domination number in non-weighted simple graphs. 
As mentioned above, for this kind of optimization problem, the probability $p_i$ of a vertex $v_i$ to be in a dominating set of a reasonably optimal weight and size should depend on the vertex weight $w_i$ and take into consideration some vertex degrees in the graph.
Let us represent this by a function $p_i=f(\bar{d},w_i)$.
Since we focus on the total weight minimization problem of a small-size dominating set $X$ in $G$, we may assume that 
the probability of a vertex $v_i$ to be in $X$ 
depends on vertex degrees $\bar{d}$ and is reciprocally proportional to the vertex weight $w_i$. 
In other words, $p_i$ can be computed by an expression of the form 
\begin{equation}
p_i=p\cdot\frac{x}{w_i},
\label{proba}
\end{equation} 
where $p$ is a coefficient depending on vertex degrees in $G$ and $x$ is a coefficient depending on vertex weights in $G$ such that
$$
\displaystyle 0\le p\cdot\frac{x}{w_i}\le 1 \;\;\hbox{  for all  }\;\; i=1,\ldots,n.
$$


In some situations, the above probabilities $p_i$ can be the same 
and not dependent on the weights. 
For example, let us assume that the weights 
do not vary too much among the vertices, 
that is, the ratio $\displaystyle w_{\rm max}/w_{\rm min}$ is reasonably close to $1$. Then, we may also assume that 
the probability $p_i$ for each vertex $v_i\in G$ to be in the dominating set $X$ does not depend on weights. 
Indeed, substituting $x=w_{\rm ave}$ in the expression $(\ref{proba})$ above, we obtain
$\displaystyle p_i=p w_{\rm ave} / w_i$.
Now, $w_i\approx w_{\rm ave}$ implies $\displaystyle p_i\approx p$ 
for every $i=1,...,n$.

In the proof of the following upper bound, which is reminiscent of the aforementioned classic bound for $\gamma(G)$, we use equal probabilities 
$p_i=p$ when applying the probabilistic method.

\begin{thm} \label{uniform}
For any graph $G$ with $\delta\ge 1$, 
$$
\gamma_w(G) \le \left(1-{\delta \over 
{(1+\delta)}^{1+1/\delta}}\right) w_G.
$$
\end{thm}

\begin{pf}
Let $A$ be a set formed by an independent choice of vertices
of $G$, where each vertex is selected with probability $p$.
Denote by $B$ the set of vertices that are not in $A$ and do not have a neighbour in $A$: $B=V(G)- N[A]$. 
Consider the set $D= A \cup B$. Clearly, $D$ is a dominating set in $G$. 
The expected weight of such a set $D$ is
\begin{eqnarray*}
\mathbb{E}[w(D)] &=& \mathbb{E}[w(A)] + \mathbb{E}[w(B)]\\
&=& \sum_{i\,=\,1}^n w_i\cdot \mathbb{P}[v_i\in A] + 
\sum_{i\,=\,1}^n w_i\cdot \mathbb{P}[v_i\in B]\\
&=&  \sum_{i\,=\,1}^n w_i\cdot p + \sum_{i\,=\,1}^n w_i\cdot(1-p)^{d_i+1}\\
&=&  \sum_{i\,=\,1}^n w_i\cdot\left(p+ (1-p)^{d_i+1}\right)\\
&\le& \sum_{i\,=\,1}^n w_i\cdot\left(p+ (1-p)^{\delta+1}\right)\\
&=&  w_G \cdot \left(p+ (1-p)^{\delta+1}\right)
\end{eqnarray*}

Minimizing the function $\psi(p) = w_G \cdot \left(p+ (1-p)^{\delta+1}\right)$, we obtain
$$
p=1-\frac{1}{(\delta+1)^{1/\delta}}.
$$
Therefore, \vspace{-0.2cm}
\begin{eqnarray*}
\mathbb{E}[w(D)] &\le& \left(1-{\delta \over {(\delta+1)}^{1+1/\delta}}\right) w_G.
\end{eqnarray*}
Since the  expectation is an average value, there exists a particular dominating set satisfying the bound, as required. 
\qed
\end{pf}

An analogue of Theorem \ref{arn} for $\gamma_w(G)$ in the case of weighted graphs easily follows from the proof of Theorem \ref{uniform}.
{Theorem~\ref{uniform} and Corollary~\ref{analogue} can be considered as simple generalizations of classic Theorems~\ref{arn} and \ref{Caro} for the smallest-cardinality dominating set problem in a graph, when all the vertex weights are equal to $1$.}

\begin{cor} \label{analogue}
For any graph $G$, 
$$
\gamma_w(G) \le {\ln(\delta+1) +1 \over
\delta+1} w_G.
$$
\end{cor}

Notice that the probability $p=1-\frac{1}{(\delta+1)^{1/\delta}}$ used in the probabilistic construction of Theorem \ref{uniform} is the same as in the probabilistic construction used in the proof of Theorem \ref{Caro} (e.g., see \cite{GPZ2013}). Therefore, in this particular case, the corresponding randomized heuristic described in Algorithm \ref{alg1} below tends to obtain small-size dominating sets satisfying the bound of Theorem \ref{Caro} at the same time.

Suppose now that $x=w_{\rm max}$ in the expression $(\ref{proba})$:
$$\displaystyle p_i=p\cdot\frac{w_{\rm max}}{w_i} \quad \mbox{for} \quad
i=1,\ldots,n.
$$ 
Because  
$0\le p\cdot\frac{w_{\rm max}}{w_i}\le 1$, 
we obtain
$$
\displaystyle p \le \frac{w_i}{w_{\rm max}}
\quad \mbox{for} \quad i=1,\ldots,n.
$$
Hence,
$
 p\le w_{\rm min}/w_{\rm max},
$
which is effectively the third assumption in the following theorem.
The second assumption  $w_{\rm max}/w_{\rm ave}\le \delta +1$ is needed
to guarantee that $p_i\ge 0$.
It may be argued that the second condition is not very restrictive in many real-life networks and it can often be overcome by removing the vertices of small degrees.

Note that the bounds of Theorem \ref{non-uniform} are not necessarily better than the previous results.
However, the heuristic based on the proof technique of this theorem often produces better results than other heuristics. This will be illustrated 
in Section \ref{sec:experiments}.


\begin{thm} \label{non-uniform}
Let $G$ be a graph such that $\delta\ge 1$, $k= w_{\rm max}/w_{\rm ave}\le \delta +1$, 
and
$p=1-\left(\frac{k}{\delta+1}\right)^{1/\delta} \le w_{\rm min}/w_{\rm max}$.
Then
$$
\gamma_w(G) \,\le\, 
n p w_{\rm max} + \sum_{i\,=\,1}^n w_i \left(1-p\right)^{d_i+1}
\,\le\,
\left(1-\frac{\delta k^{1/\delta}}{(\delta+1)^{1+1/\delta}}\right)k w_G.
$$
\end{thm}

\begin{pf}
Let $A$ be a set formed by an independent choice of vertices
of $G$, where each vertex $v_i$ is selected with some probability 
$$\displaystyle p_i=p\cdot\frac{w_{\rm max}}{w_i},
\quad 0\le p_i\le 1, 
\quad i=1,2,\ldots,n.
$$
Similar to the proof of Theorem \ref{uniform}, denote by $B$ the set of vertices that are not in $A$ and do not have a neighbour in $A$, i.e. $B=V(G)- N[A]$. 
Consider the set $D= A \cup B$. Clearly, $D$ is a dominating set in $G$. The expected weight of $D$ is
\begin{eqnarray}
\mathbb{E}[w(D)] &=& \mathbb{E}[w(A)] + \mathbb{E}[w(B)]\nonumber \\
&=& \sum_{i\,=\,1}^n w_i\cdot\mathbb{P}[v_i\in A] + \sum_{i\,=\,1}^n w_i\cdot\mathbb{P}[v_i\in B]\nonumber \\
&=&  \sum_{i\,=\,1}^n w_i \left(p\cdot \frac{w_{\rm max}}{w_i}\right) + \sum_{i\,=\,1}^n w_i \!\prod_{v_j\in N[v_i]} \left(1-p\cdot \frac{w_{\rm max}}{w_j}\right)\nonumber \\
&\le& \sum_{i\,=\,1}^n p\,w_{\rm max} + \sum_{i\,=\,1}^n w_i \!\prod_{v_j\in N[v_i]} \left(1-p\cdot \frac{w_{\rm max}}{w_{\rm max}}\right)\nonumber \\
&=& npw_{\rm max} + \sum_{i\,=\,1}^n w_i \left(1-p\right)^{d_i+1}\nonumber \\
&\le& (nw_{\rm max})p + \sum_{i\,=\,1}^n w_i \left(1-p\right)^{\delta+1}\nonumber \\
&=& (nw_{\rm max})p + \left(1-p\right)^{\delta+1} w_G \label{cor_prob}
\end{eqnarray}

Now, the value of $p$ is obtained if we minimize the function $\xi(p) = (nw_{\rm max})p + \left(1-p\right)^{\delta+1} w_G$ with respect to $p$:
$$
p=1-\left(\frac{nw_{\rm max}}{(\delta+1)w_G}\right)^{1/\delta}=1-\left(\frac{nw_{\rm max}}{(\delta+1)nw_{\rm ave}}\right)^{1/\delta}=1-\left(\frac{k}{\delta+1}\right)^{1/\delta}\!.
$$
Notice that the assumption $k=\frac{w_{\rm max}}{w_{\rm ave}}\le {\delta +1}$ guarantees that $0\le p\le 1$, 
but we additionally need to assume $p\le \frac{w_{\rm min}}{w_{\rm max}}$ to guarantee $0\le p_i\le 1$ for any $i=1,...,n$.

We already know that
$$
\gamma_w(G) \le \mathbb{E}[w(D)] \le
n p w_{\rm max} + \sum_{i\,=\,1}^n w_i \left(1-p\right)^{d_i+1}.
$$
Also,
\begin{eqnarray*}
\mathbb{E}[w(D)] &\le& pnw_{\rm max} + w_G (1-p)^{\delta+1} \\
 & = & pnw_{\rm ave}\cdot\frac{w_{\rm max}}{w_{\rm ave}} + w_G  (1-p)^{\delta+1} \\
   & = & p w_G k + w_G \left(\frac{k}{\delta+1}\right)^{1+1/\delta}\\
  & = &  \left(p  + \frac{k^{1/\delta}}{(\delta+1)^{1+1/\delta}}\right) k w_G \\
  & = & \left(1-\left(\frac{k}{\delta+1}\right)^{1/\delta}\!\! + \frac{k^{1/\delta}}{(\delta+1)^{1+1/\delta}}\right) k w_G\\
  & = & \left(1-\frac{\delta k^{1/\delta}}{(\delta+1)^{1+1/\delta}}\right)k w_G,
\end{eqnarray*}
as required. \qed
\end{pf}

Notice that we use different assumptions in the proofs of Theorems \ref{non-uniform} and \ref{uniform}. 
In the proof of Theorem \ref{uniform},
the differences in the vertex weights are ignored when trying to select vertices for the initial set $A$. 
In contrast, in the proof of Theorem \ref{non-uniform}, the probability of a vertex to be in the initial set $A$ is assumed to be reciprocally proportional to its weight, which better reflects the situation
 with the distribution of vertex weights. 
 This and some other points are illustrated later in the computational experiments in Section \ref{sec:experiments}. 

In the proof of our next theorem, we assume that 
the probability of a vertex $v_i$ to be in 
the initial set $A$ inversely depends on its weight $w_i$, $i=1,2,...,n$.
More precisely, the probability is supposed to be computed by an expression of the form 
\begin{equation}
p_i=p\cdot\left(1-\frac{w_i}{\alpha}\right),
\label{proba-inverse}
\end{equation} 
where $\alpha$ is a coefficient depending on some weights in $G$ and 
$$
0\le p\left(1-\frac{w_i}{\alpha}\right)\le 1
\quad \mbox{for all} \quad i=1,\ldots,n.
$$
If we set $\alpha=w_{\rm min}+w_{\rm max}$,
the next statement follows.


\begin{thm} \label{thm-inverse}
Let $G$ be a graph such that $\delta\ge 1$, $z= w_{\rm max}/w_{\rm min}\le \delta +1$, and $q=  1-\left({z\over \delta+1}\right)^{1/\delta}$.
Then
$$
\gamma_w(G) \,\le\, q z w_{G} 
+ \sum_{i\,=\,1}^n w_i \left(1-q\right)^{d_i+1}
\,\le\,
\left(1-\frac{\delta z^{1/\delta}}{(\delta+1)^{1+1/\delta}}\right)z w_G.
$$
\end{thm}

\begin{pf} 
Let us denote by $A$ a set formed by an independent choice of vertices
of $G$, where each vertex $v_i$ is selected with some probability 
$$\displaystyle p_i=p\left(1-\frac{w_i}{\alpha}\right), \quad 0\le p_i\le 1, \quad i=1,2,\ldots,n.$$
Similar to the proofs of Theorem \ref{uniform},
we construct the set $B=V(G)- N[A]$ and the dominating set $D= A \cup B$. 
Taking into account that
$$
\frac{w_{\rm min}}{\alpha} \le   
1-\frac{w_i}{\alpha}
\le \frac{w_{\rm max}}{\alpha},
$$
we obtain for the expected weight of $D$:
\begin{eqnarray}
\mathbb{E}[w(D)] &=& \mathbb{E}[w(A)] + \mathbb{E}[w(B)]\nonumber \\
&=& \sum_{i\,=\,1}^n w_i\cdot\mathbb{P}[v_i\in A] + \sum_{i\,=\,1}^n w_i\cdot\mathbb{P}[v_i\in B]\nonumber \\
&=&  \sum_{i\,=\,1}^n w_{i\,} p\left(1-\frac{w_i}{\alpha}\right) + \sum_{i\,=\,1}^n w_i \!\prod_{v_j\in N[v_i]} \left(1- p\left(1-\frac{w_j}{\alpha}\right)\right)\nonumber \\
&\le&  \sum_{i\,=\,1}^n w_{i\,} p \cdot \frac{w_{\rm max}}{\alpha} + \sum_{i\,=\,1}^n w_i \!\prod_{v_j\in N[v_i]} \left(1- p\cdot\frac{w_{\rm min}}{\alpha}\right)\nonumber \\
&=& \sum_{i\,=\,1}^n w_{i\,} p \cdot\frac{w_{\rm max}}{\alpha} 
+ \sum_{i\,=\,1}^n w_i \left(1-p\cdot\frac{w_{\rm min}}{\alpha}\right)^{d_i+1}\nonumber \\
&=& \frac{p z w_{G}}{z+1} 
+ \sum_{i\,=\,1}^n w_i \left(1-\frac{p}{z+1}\right)^{d_i+1}\!\!.\nonumber 
\end{eqnarray}
If we denote $q=p/(z+1)$, then
\begin{eqnarray}
\mathbb{E}[w(D)] &\le& q z w_{G} 
+ \sum_{i\,=\,1}^n w_i \left(1-q\right)^{d_i+1}\nonumber \\
&\le& q z w_{G} 
+ \sum_{i\,=\,1}^n w_i \left(1-q\right)^{\delta+1}\nonumber \\
&=& \left(qz+\left(1-q\right)^{\delta+1}\right) w_G.\nonumber 
\end{eqnarray}
It is easy to see that the last expression is minimized at
$$
q=1-\left({z\over \delta+1}\right)^{1/\delta}\!\!.
$$

Notice that to guarantee that $0\le p_i\le 1$ for all $i=1,2,...,n$, we need to require
$q=  1-\left({z\over \delta+1}\right)^{1/\delta}  \le w_{\rm min}/w_{\rm max} = \frac{1}{z}$ 
{($q\ge 0$ is guaranteed by $z\le \delta +1$)}.
However, it is possible to see that the inequality {$1-\left({z\over \delta+1}\right)^{1/\delta}  \le \frac{1}{z}$} is equivalent to 
$$\left(1-\frac{1}{z}\right)^\delta \le \frac{z}{\delta +1},$$
which holds for any $z$, $1\le z \le \delta+1$, $\delta \ge 1$. 
This is because the function $f(z)=\frac{z^{1+1/\delta}}{(\delta+1)^{1/\delta}} -z +1$ is non-negative for all $z\in [1,\delta +1]$ (it is non-negative at each endpoint of the interval and at the critical point on the interval).
Then, since naturally $z=w_{\rm max}/w_{\rm min} \ge 1$, the condition $z \le \delta+1$ is enough to have 
$0\le p_i\le 1$ for all $i=1,...,n$.
Thus, we have 
\begin{eqnarray*}
\mathbb{E}[w(D)] &\le& \left(1-\frac{\delta z^{1/\delta}}{(\delta+1)^{1+1/\delta}}\right) z w_G,
\end{eqnarray*}
as required. \qed
\end{pf}

The second and third conditions for vertex weights in Theorem \ref{non-uniform} can be rewritten as
$$\left(1-\frac{w_{\rm min}}{w_{\rm max}}\right)^\delta \cdot (\delta +1) \le k= \frac{w_{\rm max}}{w_{\rm ave}} \le \delta +1,$$
and the corresponding vertex weights conditions in Theorem \ref{thm-inverse} are 
$$1\le z=\frac{w_{\rm max}}{w_{\rm min}} \le \delta+1.$$
Clearly, there are problem instances {where} the conditions of Theorem \ref{non-uniform} are satisfied, but the conditions of Theorem \ref{thm-inverse} are not satisfied, e.g. when $k=\frac{w_{\rm max}}{w_{\rm ave}}$ is close to $\delta +1$, but $z=\frac{w_{\rm max}}{w_{\rm min}} > \delta+1$. 
On the other hand, some problem instances satisfy the conditions of Theorem \ref{thm-inverse}, but not the conditions of Theorem \ref{non-uniform}, e.g. when $z=\frac{w_{\rm max}}{w_{\rm min}}$ and $k=\frac{w_{\rm max}}{w_{\rm ave}}$ are close enough to $1$, but $k= \frac{w_{\rm max}}{w_{\rm ave}} < \left(1-\frac{w_{\rm min}}{w_{\rm max}}\right)^\delta \cdot (\delta +1)$.
More precisely, a graph with $\delta=9$, a sufficiently large number of vertices, and vertex weights distributed uniformly from $w_{\rm min}=1$ to $w_{\rm max}=10$ will have $w_{\rm ave}=5.5$ and will satisfy the conditions of Theorem \ref{thm-inverse}, but not those of Theorem \ref{non-uniform}.
 Thus, both theorems and the two corresponding probabilistic constructions are meaningful.

Theorems \ref{uniform}, \ref{non-uniform}, and \ref{thm-inverse} provide respectively probabilities\\ 
$$p_i=1-\frac{1}{(\delta+1)^{1/\delta}}\ ,$$
$$p_i=\left(1-\left(\frac{w_{\rm max}}{(\delta+1)w_{\rm ave}}\right)^{1/\delta}\right)\cdot\frac{w_{\rm max}}{w_i}\ ,\ \ \ {\rm and}$$
$$p_i=\left(1-\left({{w_{\rm max}}\over (\delta+1){w_{\rm min}}}\right)^{1/\delta}\right)\cdot \left(1+{{w_{\rm max} - w_i}\over {w_{\rm min}}}\right),\ \ \ i=1,2,\ldots,n,$$
for randomized algorithms to find reasonably good solutions (by weight) in a reasonable amount of time for large-scale instances of the problems. 
Algorithm $\ref{alg1}$ presented below is a framework for using these probabilities. 
This framework creates a  randomized algorithm from each of the theorems.
Notice that the probability $p$ used in the proof of Theorem \ref{uniform} is the same as in the probabilistic construction for the proof of Theorem \ref{Caro} and the corresponding randomized algorithm. 
This provides certain optimality in finding small-size dominating sets. 

\restylealgo{ruled}
\begin{algorithm}[h!]\label{alg1}
    \caption{Randomized small-weight dominating set}

    \KwIn{A vertex-weighted graph $G$.}
    \KwOut{A small-weight dominating set $D'$ of $G$.}
    \BlankLine
\Begin{

Compute the probability 
$\displaystyle{p_i}$ for each vertex $v_i\in V$, $i=1,\ldots,n$\;
\SetLine
Initialize set $A=\emptyset$;\\ 
 \ForEach{vertex $v_i\in V(G)$} {
    with probability $p_i$, decide whether $v_i\in A$ or $v_i\not\in A$;\\ 
    \tcc*[h]{this forms a subset $A\subseteq V(G)$}
}
\SetLine
{Use a greedy heuristic to add vertices into the set $A$ to obtain a dominating set $D$ in $G$;\\
Find a minimal by inclusion dominating set $D'\subseteq D$ in $G$;}\\
\Return $D'$}
\end{algorithm}

Computational experiments to illustrate and evaluate different solution methods are described in the next section.
The deterministic method, which uses the ILP formulations and generic ILP solvers, allows us to solve problem instances for graphs of only a few hundred vertices. Therefore, simple {and quick} heuristic solution methods provided by the probabilistic constructions in this section become important tools to tackle the problems {for larger graphs}.
The corresponding three randomized algorithms are experimentally tested in Section \ref{sec:experiments}.
Clearly, some heuristic enhancements can be used to make Algorithm $\ref{alg1}$ more effective and efficient, in particular, when forming the {initial dominating set $D$ in $G$} and finding the minimal (by inclusion) {dominating} set $D'$.


\section{Experimental evaluation}
\label{sec:experiments}

As an illustration for the results of Sections \ref{sec:deterministic-heuristic} and \ref{Section2}, 
we have implemented and tested the deterministic and heuristic solution methods 
for both problems 
on random graph instances of two types. The first type of random graphs is obtained by using the classic Erd\H{o}s--R\'enyi  random graph model \cite{G1959,ER1960}, and the other type corresponds to the random graph model used to prove asymptotic sharpness of the upper bounds of Theorems \ref{arn} and \ref{Caro} (e.g., see \cite{ZP2011}).

{For these experiments, 
in Algorithm~\ref{alg1}, after the initial set $A$ is randomly generated, we use a greedy heuristic to extend the set $A$ recursively and to obtain a dominating set $D$ of $G$.
The greedy heuristic is based on choosing a vertex adjacent to the largest number of vertices which are currently not dominated. This is similar to the greedy strategy described in \cite{CG2021}. 
A minimal by inclusion dominating set in Algorithm~\ref{alg1} is found by using the greedy heuristic described in Algortihm~\ref{MDS}. 
Notice that vertex weights are not used in the greedy strategies of both the recursive vertex selection and removal. 
From our experiments, simply running randomized Algorithm~\ref{alg1} for more iterations usually provides better results; that is, usage of more CPU time can be considered as another natural heuristic improvement for Algorithm~\ref{alg1}, e.g., if required in applications.\\}

\begin{algorithm}[H]
    \label{MDS}
    \caption{Minimal Dominating Subset}
    
    \KwIn{A graph $G=(V,E)$, a dominating set $D$ of $G$.}
    \KwOut{A minimal by inclusion dominating set $D'\subseteq D$ of $G$.}
    \BlankLine
    
    \Begin{
        Order the vertices $v_i \in D$, $i=1,\ldots,|D|$, to have $i \leq j \iff |N(v_i) - D| \leq |N(v_j) - D|$ for $1\leq i\leq j\leq |D|$;\\
        Put $D':=D$;\\
        \ForEach{vertex $v_i \in D$, $i=1,\ldots,|D|$, } {
            \If{$D' - v_i$ is a dominating set of $G$} {
                Put $D':=D' - v_i$;
            }
        }
        \Return $D'$
    }
\end{algorithm}

All the algorithms and solution methods were implemented using computer programming language C/C++, and the experiments were conducted on 
a Stone desktop PC with a $3.00$ GHz Intel Core i5 processor and $16$ GB of RAM, running Windows 10 Education OS, version 21H2. We used Gurobi Optimizer \cite{Gurobi} to solve deterministically and heuristically 
the two problems by considering their ILP formulations described in Section \ref{sec:deterministic-heuristic}.

The ILP reduction (\ref{ILP2}) to the single-objective optimization problem of finding $\gamma_{w'}(G)$ from Section \ref{sec:deterministic-heuristic} has not been used to search for heuristic randomized solutions for the two-objective optimization problem of finding $\gamma^*_w(G)$ because in the reduced problem of finding $\gamma_{w'}(G)$ the average vertex weight is $w'_{\rm ave}=1+1/n$, where $n$ is the order of the graph {(in general, $w'_i = 1+ \frac{w_i}{w_G}$, where $w_G=\sum_{i=1}^{n}w_i$, so that $1<w'_i<2$, $i=1,...,n$)}. Assuming the initial weights $w_i$ are distributed uniformly, this makes all the weights $w'_i$ in the reduced single-objective optimization problem of finding $\gamma_{w'}(G)$ very close to $1$ and not varying much among the vertices, $i=1,...,n$. Thus, the assumptions described before Theorem \ref{uniform} {-- the weights do not vary too much among the vertices, and the probabilities do not depend on weights --} are satisfied in this case. {Therefore,} the corresponding randomized algorithm can be applied to the graph with the initial vertex weights directly. On the other hand, it is possible to see that the bounds and probabilities in Theorems \ref{non-uniform} and \ref{thm-inverse} are close to those in Theorem \ref{uniform} in the case of reduction (\ref{ILP2}) and when searching for $\gamma_{w'}(G)$.

In Sections~\ref{ER-random} and \ref{sun-random}, we show that the randomized heuristics arising from the proofs of Theorems~\ref{non-uniform} and \ref{thm-inverse} are more sensitive to vertex weights in the graphs and provide better results in comparison to the randomized heuristic of Theorems~\ref{uniform}. 
Two sets of graphs are used for this comparison.
In Sections~\ref{ER-generic} and \ref{sun-generic}, we use Gurobi \cite{Gurobi} as a benchmark ILP generic solver to obtain some deterministic results for small size graphs and heuristic results for medium size graphs to show that the new randomized heuristics can find better quality initial solutions to the problems, and can do it much quicker than Gurobi. Also, these experiments show that the quality of quickly found randomized heuristic solutions is reasonably close to the Gurobi heuristic solutions, although the ILP generic solver is run on the test instances for a long time ($30$ minutes of CPU time). Finally, in Sections~\ref{ER-large} and \ref{sun-large}, we consider large-scale graphs and show that the new randomized heuristics clearly provide better results than the ILP generic solver Gurobi \cite{Gurobi}, while using less memory and CPU time resources.


\subsection{Erd\H{o}s--R\'enyi random graph model}
\label{ER-general}

For the Erd\H{o}s--R\'enyi model, denoted by $ER(n,p)$, a graph $G$ of order $n\in \mathbb{N}$ is generated in such a way that, for each (unordered) pair of distinct vertices $u$ and $v$ of $G$, the corresponding edge has the same probability $p\in [0,1]$ to be included in the graph. Such a model generates a graph whose degree distribution is a Binomial distribution with $n-1$ trials and probability $p$. In other words, given a graph order $n\in \mathbb{N}$, $n\ge 2$, and an edge-inclusion probability $p\in [0,1]$, one starts with the empty graph on $n$ vertices. Then an edge $uv$ is added to the graph with probability $p$ (or, alternatively, not added with probability $1-p$) independently for each pair of distinct vertices $\{u,v\}$. This results in a graph $G\in ER(n,p)$.

As a part of the experiments, Erd\H{o}s--R\'enyi random (unweighted) graphs on $100k$ vertices were generated for $k=1,...,10$, using the edge-inclusion probability $p=1/3$, which was manually determined to be the most illustrative in an ad-hoc way. Ten graphs were generated for each $k$, $k=1,...,10$, one hundred graphs in total, {to form the set $ER(n=100k,p=1/3)$ of $100$ test graphs}.
Then, for each graph $G\in ER(n=100k,p=1/3)$, $k=1,...,10$, weights were assigned to its vertices as integer numbers in the range $\{101,102,..., 200\}$, uniformly at random. The choice of vertex weights from the range $\{101,102,..., 200\}$ was motivated by the assumption that opening a facility must have a certain minimum basic cost ($100$ in this case) plus a certain percentage of potential additional costs. On the other hand, this range was chosen to increase the likelihood that $G$ satisfies the assumptions of Theorem \ref{non-uniform} in a simple way.

\subsubsection{{Randomized heuristics and Erd\H{o}s--R\'enyi graphs}}
\label{ER-random}

Three randomized heuristics based on Algorithm~\ref{alg1} were run on each of the above $100$ graph instances
to {quickly}  find a reasonably good solution for the minimum weight dominating set problem in $G$.
The three heuristics arise from the probabilistic constructions of Theorems \ref{uniform}, \ref{non-uniform}, and \ref{thm-inverse}, respectively. 
The aggregated averages for ten graph instances of each order $n=100k$, $k=1,...,10$, are presented in {Table \ref{table-ER-random}, 
together with the corresponding CPU run-times.}\\

\begin{table}[h!]
\centering \footnotesize 
    \begin{tabular}{|c||ccc||ccc||ccc|}
    \hline 
     \multirow{2}{*}{\begin{tabular}[c]{@{}cl@{}}$|G|$,\\ $n=100k$\end{tabular}} & \multicolumn{3}{l|}{Theorem \ref{uniform}}                                                      & \multicolumn{3}{l|}{Theorem \ref{non-uniform}}                                                      & \multicolumn{3}{l|}{Theorem \ref{thm-inverse}}                                                      \\ \cline{2-10} 
                                                                                   & \multicolumn{1}{l|}{Size} & \multicolumn{1}{l|}{Wt} & Time(s)                   & \multicolumn{1}{l|}{Size} & \multicolumn{1}{l|}{Wt} & Time(s)                   & \multicolumn{1}{l|}{Size} & \multicolumn{1}{l|}{Wt} & Time(s)                   \\ \hline \hline 
    \multirow{2}{*}{$k=1$}                                                           & \multicolumn{1}{l|}{\cellcolor{orange!25}6}    & \multicolumn{1}{l|}{897}    & \multirow{2}{*}{0.027} & \multicolumn{1}{l|}{6.1}  & \multicolumn{1}{l|}{825.5}  & \multirow{2}{*}{0.027} & \multicolumn{1}{l|}{\cellcolor{orange!25}6}    & \multicolumn{1}{l|}{835.6}  & \multirow{2}{*}{0.027} \\ \cline{2-3} \cline{5-6} \cline{8-9}
                                                                                   & \multicolumn{1}{l|}{6.3}  & \multicolumn{1}{l|}{885.2}  &                           & \multicolumn{1}{l|}{6.2}  & \multicolumn{1}{l|}{\cellcolor{orange!25}824.4}  &                           & \multicolumn{1}{l|}{6.2}  & \multicolumn{1}{l|}{833.1}  &                           \\ \hline
    \multirow{2}{*}{$k=2$}                                                           & \multicolumn{1}{l|}{7.5}  & \multicolumn{1}{l|}{1066.9} & \multirow{2}{*}{0.062} & \multicolumn{1}{l|}{\cellcolor{orange!25}7.4}  & \multicolumn{1}{l|}{1041}   & \multirow{2}{*}{0.061}  & \multicolumn{1}{l|}{7.5}  & \multicolumn{1}{l|}{1049.8} & \multirow{2}{*}{0.062} \\ \cline{2-3} \cline{5-6} \cline{8-9}
                                                                                   & \multicolumn{1}{l|}{7.8}  & \multicolumn{1}{l|}{1047.6} &                           & \multicolumn{1}{l|}{7.6}  & \multicolumn{1}{l|}{\cellcolor{orange!25}1011.8} &                           & \multicolumn{1}{l|}{7.7}  & \multicolumn{1}{l|}{1027.8} &                           \\ \hline
    \multirow{2}{*}{$k=3$}                                                           & \multicolumn{1}{l|}{8.6}  & \multicolumn{1}{l|}{1225.5} & \multirow{2}{*}{0.1} & \multicolumn{1}{l|}{8.5}  & \multicolumn{1}{l|}{1140.4} & \multirow{2}{*}{0.1} & \multicolumn{1}{l|}{\cellcolor{orange!25}8.4}  & \multicolumn{1}{l|}{1217.7} & \multirow{2}{*}{0.1} \\ \cline{2-3} \cline{5-6} \cline{8-9}
                                                                                   & \multicolumn{1}{l|}{8.6}  & \multicolumn{1}{l|}{1225.5} &                           & \multicolumn{1}{l|}{8.7}  & \multicolumn{1}{l|}{\cellcolor{orange!25}1126.6} &                           & \multicolumn{1}{l|}{8.8}  & \multicolumn{1}{l|}{1185.1} &                           \\ \hline
    \multirow{2}{*}{$k=4$}                                                           & \multicolumn{1}{l|}{\cellcolor{orange!25}9}    & \multicolumn{1}{l|}{1293.6} & \multirow{2}{*}{0.16}  & \multicolumn{1}{l|}{9.3}  & \multicolumn{1}{l|}{1306}   & \multirow{2}{*}{0.15}  & \multicolumn{1}{l|}{9.1}  & \multicolumn{1}{l|}{1308.4} & \multirow{2}{*}{0.16}  \\ \cline{2-3} \cline{5-6} \cline{8-9}
                                                                                   & \multicolumn{1}{l|}{9.2}  & \multicolumn{1}{l|}{\cellcolor{orange!25}1276.4} &                           & \multicolumn{1}{l|}{9.5}  & \multicolumn{1}{l|}{1290.3} &                           & \multicolumn{1}{l|}{9.5}  & \multicolumn{1}{l|}{1283.1} &                           \\ \hline
    \multirow{2}{*}{$k=5$}                                                           & \multicolumn{1}{l|}{9.8}  & \multicolumn{1}{l|}{1387.5} & \multirow{2}{*}{0.21}  & \multicolumn{1}{l|}{10}   & \multicolumn{1}{l|}{1346.9} & \multirow{2}{*}{0.21}   & \multicolumn{1}{l|}{\cellcolor{orange!25}9.7}  & \multicolumn{1}{l|}{1341.8} & \multirow{2}{*}{0.21}  \\ \cline{2-3} \cline{5-6} \cline{8-9}
                                                                                   & \multicolumn{1}{l|}{9.9}  & \multicolumn{1}{l|}{1381.4} &                           & \multicolumn{1}{l|}{10.1} & \multicolumn{1}{l|}{1345.8} &                           & \multicolumn{1}{l|}{9.8}  & \multicolumn{1}{l|}{\cellcolor{orange!25}1340.5} &                           \\ \hline
    \multirow{2}{*}{$k=6$}                                                           & \multicolumn{1}{l|}{\cellcolor{orange!25}10}   & \multicolumn{1}{l|}{1433.5} & \multirow{2}{*}{0.29}  & \multicolumn{1}{l|}{10.1} & \multicolumn{1}{l|}{1385.7} & \multirow{2}{*}{0.28}  & \multicolumn{1}{l|}{\cellcolor{orange!25}10}   & \multicolumn{1}{l|}{1403.1} & \multirow{2}{*}{0.28}  \\ \cline{2-3} \cline{5-6} \cline{8-9}
                                                                                   & \multicolumn{1}{l|}{10.2} & \multicolumn{1}{l|}{1431.3} &                           & \multicolumn{1}{l|}{10.2} & \multicolumn{1}{l|}{\cellcolor{orange!25}1384.1} &                           & \multicolumn{1}{l|}{10.5} & \multicolumn{1}{l|}{1392.6} &                           \\ \hline
    \multirow{2}{*}{$k=7$}                                                           & \multicolumn{1}{l|}{10.5} & \multicolumn{1}{l|}{1463.9} & \multirow{2}{*}{0.37}  & \multicolumn{1}{l|}{\cellcolor{orange!25}10.4} & \multicolumn{1}{l|}{1430.7} & \multirow{2}{*}{0.36}    & \multicolumn{1}{l|}{10.7} & \multicolumn{1}{l|}{1462.8} & \multirow{2}{*}{0.36}  \\ \cline{2-3} \cline{5-6} \cline{8-9}
                                                                                   & \multicolumn{1}{l|}{10.5} & \multicolumn{1}{l|}{1463.9} &                           & \multicolumn{1}{l|}{10.4} & \multicolumn{1}{l|}{\cellcolor{orange!25}1430.7} &                           & \multicolumn{1}{l|}{10.7} & \multicolumn{1}{l|}{1462.8} &                           \\ \hline
    \multirow{2}{*}{$k=8$}                                                           & \multicolumn{1}{l|}{11}   & \multicolumn{1}{l|}{1599.7} & \multirow{2}{*}{0.45}  & \multicolumn{1}{l|}{11}   & \multicolumn{1}{l|}{1532.7} & \multirow{2}{*}{0.45}  & \multicolumn{1}{l|}{\cellcolor{orange!25}10.8} & \multicolumn{1}{l|}{1538.9} & \multirow{2}{*}{0.45}  \\ \cline{2-3} \cline{5-6} \cline{8-9}
                                                                                   & \multicolumn{1}{l|}{11.1} & \multicolumn{1}{l|}{1599.1} &                           & \multicolumn{1}{l|}{11.3} & \multicolumn{1}{l|}{1524.2} &                           & \multicolumn{1}{l|}{11.1} & \multicolumn{1}{l|}{\cellcolor{orange!25}1507.5} &                           \\ \hline
    \multirow{2}{*}{$k=9$}                                                           & \multicolumn{1}{l|}{11.1} & \multicolumn{1}{l|}{1626.3} & \multirow{2}{*}{0.56}   & \multicolumn{1}{l|}{11.1} & \multicolumn{1}{l|}{1591.6} & \multirow{2}{*}{0.56}  & \multicolumn{1}{l|}{\cellcolor{orange!25}11}   & \multicolumn{1}{l|}{1557}   & \multirow{2}{*}{0.54}  \\ \cline{2-3} \cline{5-6} \cline{8-9}
                                                                                   & \multicolumn{1}{l|}{11.5} & \multicolumn{1}{l|}{1605.1} &                           & \multicolumn{1}{l|}{11.6} & \multicolumn{1}{l|}{1556.2} &                           & \multicolumn{1}{l|}{11.2} & \multicolumn{1}{l|}{\cellcolor{orange!25}1531.2} &                           \\ \hline
    \multirow{2}{*}{$k=10$}                                                          & \multicolumn{1}{l|}{\cellcolor{orange!25}11.3} & \multicolumn{1}{l|}{1622.6} & \multirow{2}{*}{0.64}  & \multicolumn{1}{l|}{11.4} & \multicolumn{1}{l|}{1542}   & \multirow{2}{*}{0.64}  & \multicolumn{1}{l|}{\cellcolor{orange!25}11.3} & \multicolumn{1}{l|}{1590.2} & \multirow{2}{*}{0.64}  \\ \cline{2-3} \cline{5-6} \cline{8-9}
                                                                                   & \multicolumn{1}{l|}{11.4} & \multicolumn{1}{l|}{1622}   &                           & \multicolumn{1}{l|}{11.6} & \multicolumn{1}{l|}{\cellcolor{orange!25}1534.6} &                           & \multicolumn{1}{l|}{11.6} & \multicolumn{1}{l|}{1573.8} &                           \\ \hline
    \end{tabular}
    \caption{\label{table-ER-random}Aggregated results of running the randomized heuristics on the Erd\H{o}s--R\'enyi graphs.} 
\end{table}

Each of the randomized algorithms was run {twenty} times on each graph instance, and the best found solution (out of {twenty}) by the dominating set size and also by the dominating set weight were recorded. 
The corresponding average solution parameters for the best found set size are shown in the upper subrows, and the averages for the best found solutions by the set weight are shown in the lower subrows, 
for each $k=1,...,10$, in {Table \ref{table-ER-random}}. 
Each of the randomized algorithms of Theorems \ref{uniform}, \ref{non-uniform}, and \ref{thm-inverse} finds better solutions by the set size for {$73\%, 68\%$, and $77\%$} of all instances, respectively (corresponding dominating sets may have the same cardinality). In other words, the three randomized algorithms show a similar performance by the dominating set size. However, the randomized algorithm arising from the proof of Theorem \ref{uniform} is much less successful in finding the best heuristic solution by the set weight (only {$20\%$} of all instances). The randomized algorithms arising from the proofs of Theorems \ref{non-uniform} and \ref{thm-inverse} find better solutions by the set weight for {$45\%$ and $36\%$} of all instances, respectively, i.e. perform clearly better for this parameter.

The upper bounds of Theorems \ref{Caro} and \ref{uniform} were satisfied for all the problem instances, with the results getting closer to the upper bounds for graphs of larger order. 
It is possible to see from {Table \ref{table-ER-random}} that, in most of the cases, the heuristic methods derived from Theorems \ref{non-uniform} and \ref{thm-inverse} provide better results for the two-objective optimization problem by weight and by size than that of Theorem \ref{uniform}, although both methods of Theorems \ref{non-uniform} and \ref{thm-inverse} are designed for optimization by weight only.

\subsubsection{{Using a generic solver on Erd\H{o}s--R\'enyi graphs}}
\label{ER-generic}

For each weighted graph $G\in ER(n=100k,p=1/3)$, $k=1,...,10$, {considered above,} we made an attempt to find exact deterministic solutions to the problems of computing $\gamma^*_w(G)$ and $\gamma_w(G)$ by using the ILP formulations 
described in Section \ref{sec:deterministic-heuristic}
and {the ILP generic solver Gurobi \cite{Gurobi}}. 
This was successful in a reasonable amount of {CPU time of at most $30$ minutes ($1800$ sec)} only for $k=1,2$, {when computing $\gamma^*_w(G)$, and for $k=1,2,3$, when computing $\gamma_w(G)$}. 
In these deterministic computational experiments, only one graph instance was found to have a dominating set of a non-minimum size with weight lower  than $\gamma^*_w(G)$, i.e. with $\gamma_w(G)<\gamma^*_w(G)$ (see Table \ref{table-ER-Gurobi}).
\noindent
\begin{table}[h!]
    \centering \footnotesize 
    \begin{tabular}{|c||ccclc||ccclc|}
    \hline
    \multirow{2}{*}{\begin{tabular}[c]{c@{}c@{}}$|G|$,\\ $n=100k$\end{tabular}} & \multicolumn{5}{c|}{ILP for $\gamma_w^*(G)$}                                                                                                                                                                                                                                                        & \multicolumn{5}{c|}{ILP for $\gamma_w(G)$}                                                                                                                                                                                                                                                       \\ \cline{2-11} 
                                                                                   & \multicolumn{1}{c|}{\begin{tabular}[c]{@{}c@{}}IHS\\ size\end{tabular}} & \multicolumn{1}{c|}{\begin{tabular}[c]{@{}c@{}}IHS CPU\\  time\,(s)\end{tabular}} & \multicolumn{1}{c|}{\begin{tabular}[c]{@{}c@{}}BPS\\ size\end{tabular}} & \multicolumn{1}{l|}{\begin{tabular}[c]{@{}c@{}}BPS\\ wt\end{tabular}} & \begin{tabular}[c]{@{}c@{}}BPS CPU\\ time\,(s)\end{tabular}   & \multicolumn{1}{c|}{\begin{tabular}[c]{@{}c@{}}IHS\\ wt\end{tabular}} & \multicolumn{1}{c|}{\begin{tabular}[c]{@{}c@{}}IHS CPU\\ time\,(s)\end{tabular}} & \multicolumn{1}{c|}{\begin{tabular}[c]{@{}c@{}}BPS\\ wt\end{tabular}} & \multicolumn{1}{l|}{\begin{tabular}[c]{@{}c@{}}BPS\\ size\end{tabular}} & \begin{tabular}[c]{@{}c@{}}BPS CPU\\ time\,(s)\end{tabular}    \\ \hline \hline
    $k=1$                                                                            & \multicolumn{1}{c|}{7.9}                                                      & \multicolumn{1}{c|}{0}                                                      & \multicolumn{1}{c|}{\cellcolor{orange!25}5}                                                      & \multicolumn{1}{l|}{\cellcolor{orange!25}620.4}  & \cellcolor{orange!25}1.33 & \multicolumn{1}{c|}{1205.5}                                                   & \multicolumn{1}{c|}{0}                                                      & \multicolumn{1}{c|}{\cellcolor{orange!25}616}                                                      & \multicolumn{1}{l|}{\cellcolor{orange!25}5.1}  & \cellcolor{orange!25}0.4  \\ \hline
    $k=2$                                                                            & \multicolumn{1}{c|}{9.6}                                                      & \multicolumn{1}{c|}{0.025}                                                  & \multicolumn{1}{c|}{\cellcolor{orange!25}6}                                                      & \multicolumn{1}{l|}{\cellcolor{orange!25}699.3}  & \cellcolor{orange!25}494.34 & \multicolumn{1}{c|}{1432.5}                                                   & \multicolumn{1}{c|}{0.022}                                                  & \multicolumn{1}{c|}{\cellcolor{orange!25}699.3}                                                    & \multicolumn{1}{l|}{\cellcolor{orange!25}6}    & \cellcolor{orange!25}30.42 \\ \hline
    $k=3$                                                                            & \multicolumn{1}{c|}{10.5}                                                     & \multicolumn{1}{c|}{0.084}                                                  & \multicolumn{1}{c|}{6.8}                                                    & \multicolumn{1}{l|}{800.4}  & 1800   & \multicolumn{1}{c|}{1600.6}                                                   & \multicolumn{1}{c|}{0.081}                                                  & \multicolumn{1}{c|}{\cellcolor{orange!25}762.1}                                                    & \multicolumn{1}{l|}{\cellcolor{orange!25}7}    & \cellcolor{orange!25}646.74  \\ \hline
    $k=4$                                                                            & \multicolumn{1}{c|}{10.7}                                                     & \multicolumn{1}{c|}{0.19}                                                  & \multicolumn{1}{c|}{7}                                                      & \multicolumn{1}{l|}{887.6}  & 1800   & \multicolumn{1}{c|}{1550.3}                                                   & \multicolumn{1}{c|}{0.2}                                                  & \multicolumn{1}{c|}{818.1}                                                    & \multicolumn{1}{l|}{7.4}  & 1800    \\ \hline
    $k=5$                                                                            & \multicolumn{1}{c|}{11.1}                                                     & \multicolumn{1}{c|}{0.35}                                                  & \multicolumn{1}{c|}{7.8}                                                    & \multicolumn{1}{l|}{933.9}  & 1800   & \multicolumn{1}{c|}{1651}                                                     & \multicolumn{1}{c|}{0.39}                                                   & \multicolumn{1}{c|}{859.7}                                                    & \multicolumn{1}{l|}{8}    & 1800    \\ \hline
    $k=6$                                                                            & \multicolumn{1}{c|}{11.8}                                                     & \multicolumn{1}{c|}{0.74}                                                  & \multicolumn{1}{c|}{8}                                                      & \multicolumn{1}{l|}{957.4}  & 1800   & \multicolumn{1}{c|}{1771.2}                                                   & \multicolumn{1}{c|}{0.91}                                                   & \multicolumn{1}{c|}{910.6}                                                    & \multicolumn{1}{l|}{8.2}  & 1800    \\ \hline
    $k=7$                                                                            & \multicolumn{1}{c|}{12}                                                       & \multicolumn{1}{c|}{1.14}                                                  & \multicolumn{1}{c|}{8}                                                      & \multicolumn{1}{l|}{1017.4} & 1800   & \multicolumn{1}{c|}{1827.6}                                                   & \multicolumn{1}{c|}{1.37}                                                  & \multicolumn{1}{c|}{931.6}                                                    & \multicolumn{1}{l|}{8.9}  & 1800    \\ \hline
    $k=8$                                                                            & \multicolumn{1}{c|}{12.6}                                                     & \multicolumn{1}{c|}{1.63}                                                  & \multicolumn{1}{c|}{8.8}                                                    & \multicolumn{1}{l|}{1053.1} & 1800   & \multicolumn{1}{c|}{1990.9}                                                   & \multicolumn{1}{c|}{1.55}                                                  & \multicolumn{1}{c|}{955.3}                                                    & \multicolumn{1}{l|}{9}    & 1800    \\ \hline
    $k=9$                                                                            & \multicolumn{1}{c|}{12.6}                                                     & \multicolumn{1}{c|}{2.26}                                                  & \multicolumn{1}{c|}{9}                                                      & \multicolumn{1}{l|}{1076.6} & 1800   & \multicolumn{1}{c|}{1851.4}                                                   & \multicolumn{1}{c|}{2.04}                                                  & \multicolumn{1}{c|}{972.6}                                                    & \multicolumn{1}{l|}{9}    & 1800    \\ \hline
    $k=10$                                                                           & \multicolumn{1}{c|}{13}                                                       & \multicolumn{1}{c|}{3.004}                                                  & \multicolumn{1}{c|}{9}                                                      & \multicolumn{1}{l|}{1127.4} & 1800   & \multicolumn{1}{c|}{1954.7}                                                   & \multicolumn{1}{c|}{2.92}                                                  & \multicolumn{1}{c|}{1003.5}                                                   & \multicolumn{1}{l|}{9.4}  & 1800    \\ \hline
    \end{tabular}
     \caption{\label{table-ER-Gurobi}Aggregated results of running the ILP generic solver on the Erd\H{o}s--R\'enyi graphs.}
\end{table}

For the larger (medium) size graph instances, i.e. for $k\ge 3$, {when computing $\gamma^*_w(G)$, and for $k\ge 4$, when computing $\gamma_w(G)$}, Gurobi was run for $30$ minutes ($1800$ sec CPU time) as a benchmark heuristic solver. This is to see how the Gurobi's initial heuristic solution (IHS) compares to the randomized algorithms solutions from the previous Section~\ref{ER-random}. Also, considering CPU time as a limited computational resource, these computations were run to obtain the best possible solution (BPS) of the ILP formulation in $1800$ sec CPU time. The aggregated results are presented in Table~\ref{table-ER-Gurobi} (the deterministic solution results are shaded).

It can be seen that, for the medium size graphs ($k\ge 4$), the new randomized heuristics clearly provide better results (see Table~\ref{table-ER-random}) using much less CPU run-time than the IHS of Gurobi (see IHS columns in Table~\ref{table-ER-Gurobi}).
It can also be seen that, on average, the best solutions found by the randomized algorithms for {$k=1,2,3$} are about {$34-48\%$} worse by weight and, {for $k=1,2$,} about {$20-23\%$} worse by size than the optimal {(deterministic)} solutions {(see shaded cells in BPS columns in Table~\ref{table-ER-Gurobi})}. 
{In comparison to the heuristic results for medium size graphs ($k\ge 4$) obtained by Gurobi in 1800 sec, the very quick (less than a second) heuristic solutions by the new randomized algorithms are about $52-58\%$ worse by weight and about $22-30\%$ worse by size. However, this gap can be significantly reduced and even better results can be obtained by running the randomized algorithms for more iterations (e.g., for 1800 sec CPU time). This is clearly demonstrated for large size graphs in the next Section~\ref{ER-large}.}
The average CPU run-times to compute exact values of $\gamma^*_w(G)$ and $\gamma_w(G)$ (for $k=1,2,3$) using the ILP formulations and {generic solver 
Gurobi (see BPS columns in Table~\ref{table-ER-Gurobi})} clearly show the considerable growth of computational time requirements with respect to the graph order, as well as differences in computational complexity for finding $\gamma^*_w(G)$ and $\gamma_w(G)$.

\subsubsection{{Large-scale Erd\H{o}s--R\'enyi graphs}}
\label{ER-large}

Finally, we have run the three randomized algorithms and the ILP generic solver Gurobi on two large size Erd\H{o}s--R\'enyi random graphs, generated as described at the beginning of  Section~\ref{ER-general}. One graph has $n=20,000$ vertices, and the other has $n=40,000$ vertices. Both graphs were generated with the edge probability $p=0.1$.
In these computational experiments, each of the four solvers was given $30$ minutes (1800 sec) of CPU time to find a solution to the problems corresponding to $\gamma_w^*(G)$ and $\gamma_w(G)$.

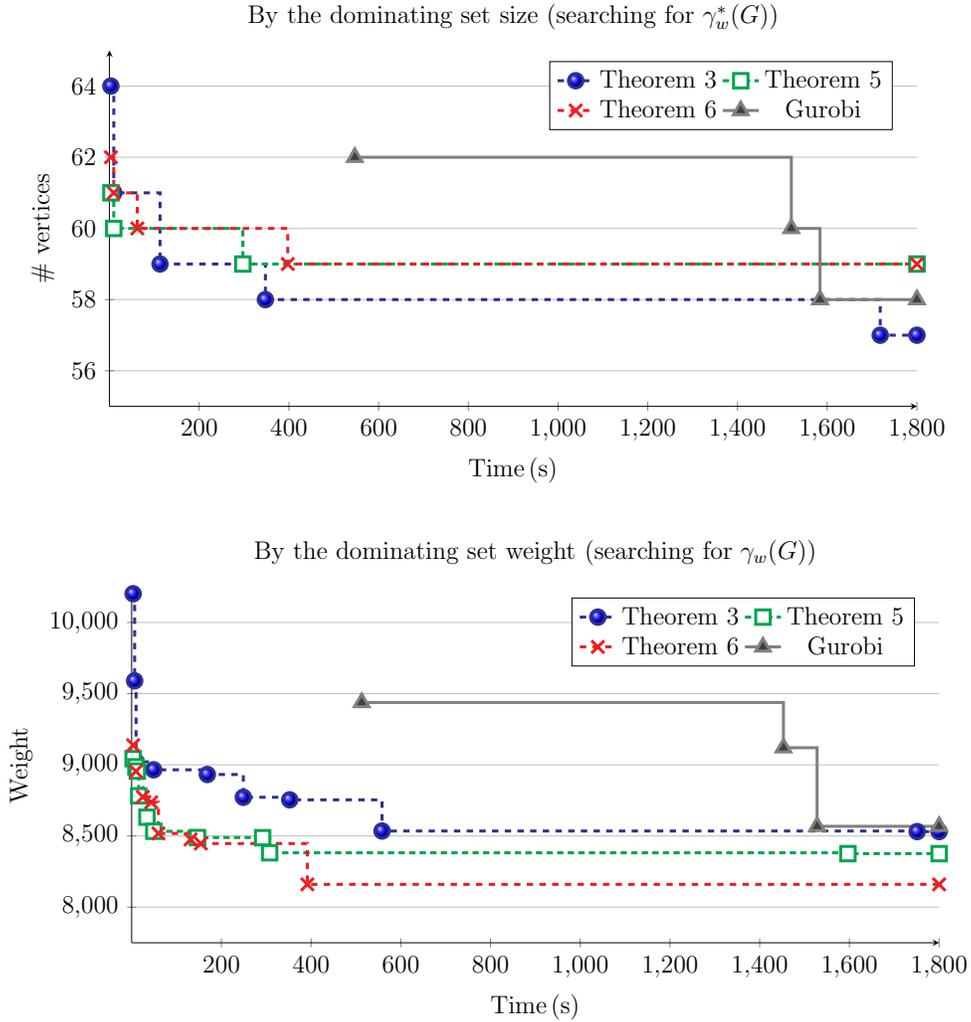
\begin{figure}[H]
    \centering
    \begin{tikzpicture}[scale=0.8]
        \begin{axis}[
            axis lines=left,
            title={By the dominating set size\  (searching for $\gamma^*_w(G)$)},
            xlabel={Time\,(s)},
            ylabel={\# vertices},
            xmin=0, xmax=1800,
            ymin=55, ymax=65,
            xtick={200,400,600,800,1000,1200,1400,1600,1800},
            ytick={56,58,60,62,64},
            ymajorgrids=true,
            legend pos=north east,
            legend columns=2,
        ]
        
        \addplot+[
            const plot,
            color=Blue,
            mark=ball,
            mark options={fill=Blue,style=solid},
            mark size=1.2mm,
            style=dashed,
            line width=0.5mm
            ] coordinates {
                (3.12547,64)(9.62558,61)(112.718,59)(347.638,58)(1718.17,57)(1800,57)
            };
            \addlegendentry{Theorem 3}
        
        \addplot+[
            const plot,
            color=Green,
            mark=square*,
            mark options={fill=White,style=solid},
            mark size=1.2mm,
            style=densely dashed,
            line width=0.5mm
            ] coordinates {
                (3.22307,61)(9.54669,60)(297.345,59)(1800,59)
            };
            \addlegendentry{Theorem 5}

        \addplot+[
            const plot,
            color=Red,
            mark=x,
            mark options={fill=Red,style=solid},
            mark size=1.5mm,
            style=dashed,
            line width=0.5mm
            ] coordinates {
                (3.28601,62)(9.48158,61)(62.0447,60)(397.443,59)(1800,59)
            };
            \addlegendentry{Theorem 6}

        \addplot+[
            const plot,
            color=gray,
            mark=triangle*,
            mark options={fill=darkgray,style=solid},
            mark size=1.4mm,
            line width=0.5mm
            ] coordinates {
                (546.85,62)(1520,60)(1584,58)(1800,58)
            };
            \addlegendentry{Gurobi}
        \end{axis}
    \end{tikzpicture}
    
\vspace{5mm}
    
    \begin{tikzpicture}[scale=0.8]
        \begin{axis}[
            axis lines=left,
            title={By the dominating set weight\ (searching for $\gamma_w(G)$)},
            xlabel={Time\,(s)},
            ylabel={Weight},
            xmin=0, xmax=1800,
            ymin=7750, ymax=10250,
            xtick={200,400,600,800,1000,1200,1400,1600,1800},
            ytick={8000,8500,9000,9500,10000},
            ymajorgrids=true,
            legend pos=north east,
            legend columns=2,
        ]
        
        \addplot+[
            const plot,
            color=Blue,
            mark=ball,
            mark options={fill=Blue,style=solid},
            mark size=1.2mm,
            style=dashed,
            line width=0.5mm
            ] coordinates {
                (3.26023,10202)(6.66307,9590)(9.90724,9021)(48.9844,8965)(168.738,8932)
                (248.965,8772)(352.42,8754)(558.161,8535)(1751.08,8531)(1800,8531)
            };
            \addlegendentry{Theorem 3}
        
        \addplot+[
            const plot,
            color=Green,
            mark=square*,
            mark options={fill=White,style=solid},
            mark size=1.2mm,
            style=densely dashed,
            line width=0.5mm
            ] coordinates {
                (3.2404,9043)(9.44757,8986)(12.6146,8956)(15.652,8782)(34.0499,8632)(49.414,8533)
                (146.784,8489)(292.054,8488)(307.73,8382)(1596.74,8376)(1800,8376)
            };
            \addlegendentry{Theorem 5}

        \addplot+[
            const plot,
            color=Red,
            mark=x,
            mark options={fill=Red,style=solid},
            mark size=1.5mm,
            style=dashed,
            line width=0.5mm
            ] coordinates {
                (3.32095,9138)(9.51018,8954)(24.9831,8775)(43.6934,8737)(59.687,8518)
                (131.43,8476)(154.214,8447)(391.65,8160)(1800,8160)
            };
            \addlegendentry{Theorem 6}

        \addplot+[
            const plot,
            color=gray,
            mark=triangle*,
            mark options={fill=darkgray,style=solid},
            mark size=1.4mm,
            line width=0.5mm
            ] coordinates {
                (512.88,9438)(1453,9120)(1528,8568)(1800,8568)
            };
            \addlegendentry{Gurobi}
        \end{axis}
    \end{tikzpicture}
    \caption{\label{fig-ER-large-S}Performance profiles of the solvers for the Erd\H{o}s--R\'enyi graph on $20,000$ vertices.}
\end{figure}

For the Erd\H{o}s--R\'enyi graph on $n=20,000$ vertices, the randomized algorithms arising from Theorems~\ref{uniform}, \ref{non-uniform}, and \ref{thm-inverse} respectively used  $566$, $558$, and $556$ iterations.
For this graph, the ILP generic solver Gurobi produced lower bounds on the solutions by size, i.e. established that $\gamma(G) \ge 10$ (in $1519$ sec CPU time), and by weight, i.e. established that $\gamma_w(G) \ge 1159$ (in 1452 sec CPU time). The performance profiles of the four tested solvers on this large graph are presented in Figure~\ref{fig-ER-large-S} and Table~\ref{table-ER-large}.


\begin{table}[h!]
    \centering \footnotesize 
    \begin{tabular}{|c||ccc|ccc|ccc|ccc|}
    \hline
    \multirow{2}{*}{Problem} & \multicolumn{3}{c|}{Theorem~\ref{uniform}}                                                                                 & \multicolumn{3}{c|}{Theorem~\ref{non-uniform}}                                                                                 & \multicolumn{3}{c|}{Theorem~\ref{thm-inverse}}                                                                                 & \multicolumn{3}{c|}{Gurobi}                                                                                    \\ \cline{2-13} 
                           & \multicolumn{1}{c|}{Size} & \multicolumn{1}{c|}{Wt} & \begin{tabular}[c]{@{}c@{}}Time\\ (sec)\end{tabular} & \multicolumn{1}{c|}{Size} & \multicolumn{1}{c|}{Wt} & \begin{tabular}[c]{@{}c@{}}Time\\ (sec)\end{tabular} & \multicolumn{1}{c|}{Size} & \multicolumn{1}{c|}{Wt} & \begin{tabular}[c]{@{}c@{}}Time\\ (sec)\end{tabular} & \multicolumn{1}{c|}{Size} & \multicolumn{1}{c|}{Wt} & \begin{tabular}[c]{@{}c@{}}Time\\ (sec)\end{tabular} \\ \hline \hline
    $\gamma^*_w(G)$                   & \multicolumn{1}{c|}{57}   & \multicolumn{1}{c|}{8531}   & 1718                                                 & \multicolumn{1}{c|}{59}   & \multicolumn{1}{c|}{8382}   & 297                                                  & \multicolumn{1}{c|}{59}   & \multicolumn{1}{c|}{8160}   & 392                                                  & \multicolumn{1}{c|}{58}   & \multicolumn{1}{c|}{8183}   & 1584                                                 \\ \hline
   $\gamma_w(G)$             & \multicolumn{1}{c|}{57}   & \multicolumn{1}{c|}{8531}   & 1718                                                 & \multicolumn{1}{c|}{60}   & \multicolumn{1}{c|}{8376}   & 1597                                                 & \multicolumn{1}{c|}{59}   & \multicolumn{1}{c|}{8160}   & 392                                                  & \multicolumn{1}{c|}{58}   & \multicolumn{1}{c|}{8568}   & 1528                                                 \\ \hline
    \end{tabular}
    \caption{\label{table-ER-large}Best found results for the Erd\H{o}s--R\'enyi graph on $20,000$ veritces.}
\end{table}


In general, the randomized algorithms find better quality solutions and much quicker than Gurobi. When trying to find $\gamma^*_w(G)$, the randomized algorithm arising from Theorem~\ref{uniform} provides the best solution by the dominating set size. This can be explained by the same optimal probability used in Theorems~\ref{uniform} and \ref{Caro}. 
When searching for $\gamma_w(G)$, all three randomized algorithms find better solutions and quicker than Gurobi. In this case, the randomized algorithm arising from Theorem~\ref{thm-inverse} provides the best solution, which can be explained by its sensitivity to the vertex weights in the graph. Notice that Gurobi finds a better solution by weight when searching for $\gamma^*_w(G)$, which is still not the best (out of four).

\begin{figure}[H]
    \centering
    \begin{tikzpicture}[scale=0.8]
        \begin{axis}[
            axis lines=left,
            title={By the dominating set size\  (searching for $\gamma^*_w(G)$)},
            xlabel={Time\,(s)},
            ylabel={\# vertices},
            xmin=0, xmax=1800,
            ymin=63, ymax=73,
            xtick={200,400,600,800,1000,1200,1400,1600,1800},
            ytick={64,66,68,70,72},
            ymajorgrids=true,
            legend pos=north east,
        ]
        
        \addplot+[
            const plot,
            color=Blue,
            mark=ball,
            mark options={fill=Blue,style=solid},
            mark size=1.2mm,
            style=dashed,
            line width=0.5mm
            ] coordinates {
                (16.0048,72)(31.2195,69)(136.81,67)(240.631,66)(1730.02,65)(1800,65)
            };
            \addlegendentry{Theorem 3}
        
        \addplot+[
            const plot,
            color=Green,
            mark=square*,
            mark options={fill=White,style=solid},
            mark size=1.2mm,
            style=densely dashed,
            line width=0.5mm
            ] coordinates {
                (15.5314,68)(124.616,67)(266.075,66)(1616.04,65)(1800,65)
            };
            \addlegendentry{Theorem 5}

        \addplot+[
            const plot,
            color=Red,
            mark=x,
            mark options={fill=Red,style=solid},
            mark size=1.5mm,
            style=dashed,
            line width=0.5mm
            ] coordinates {
                (14.7563,69)(31.1729,68)(323.962,67)(1117.87,65)(1800,65)
            };
            \addlegendentry{Theorem 6}
        \end{axis}
    \end{tikzpicture}

\vspace{5mm}

    \begin{tikzpicture}[scale=0.8]
        \begin{axis}[
            axis lines=left,
            title={By the dominating set weight\  (searching for $\gamma_w(G)$)},
            xlabel={Time\,(s)},
            ylabel={Weight},
            xmin=0, xmax=1800,
            ymin=8750, ymax=11250,
            xtick={200,400,600,800,1000,1200,1400,1600,1800},
            ytick={9000,9500,10000,10500,11000},
            ymajorgrids=true,
            legend pos=north east,
        ]
        
        \addplot+[
            const plot,
            color=Blue,
            mark=ball,
            mark options={fill=Blue,style=solid},
            mark size=1.2mm,
            style=dashed,
            line width=0.5mm
            ] coordinates {
                (16.0048,11158)(31.2195,10391)(107.553,10244)(136.81,10005)(270.1,9707)(1730.02,9284)(1800,9284)
            };
            \addlegendentry{Theorem 3}
        
        \addplot+[
            const plot,
            color=Green,
            mark=square*,
            mark options={fill=White,style=solid},
            mark size=1.2mm,
            style=densely dashed,
            line width=0.5mm
            ] coordinates {
                (15.5314,10276)(30.7991,9524)(266.075,9430)(1616.04,9418)(1800,9418)
            };
            \addlegendentry{Theorem 5}

        \addplot+[
            const plot,
            color=Red,
            mark=x,
            mark options={fill=Red,style=solid},
            mark size=1.5mm,
            style=dashed,
            line width=0.5mm
            ] coordinates {
                (14.7563,10570)(31.1729,9641)(118.403,9433)(323.962,9400)(1316.6,9386)(1458.7,9207)(1800,9207)
            };
            \addlegendentry{Theorem 6}
        \end{axis}
    \end{tikzpicture}
        \caption{\label{fig-ER-large-40-S}Running the randomized heuristics on the Erd\H{o}s--R\'enyi graph of $40,000$ vertices.}
\end{figure}
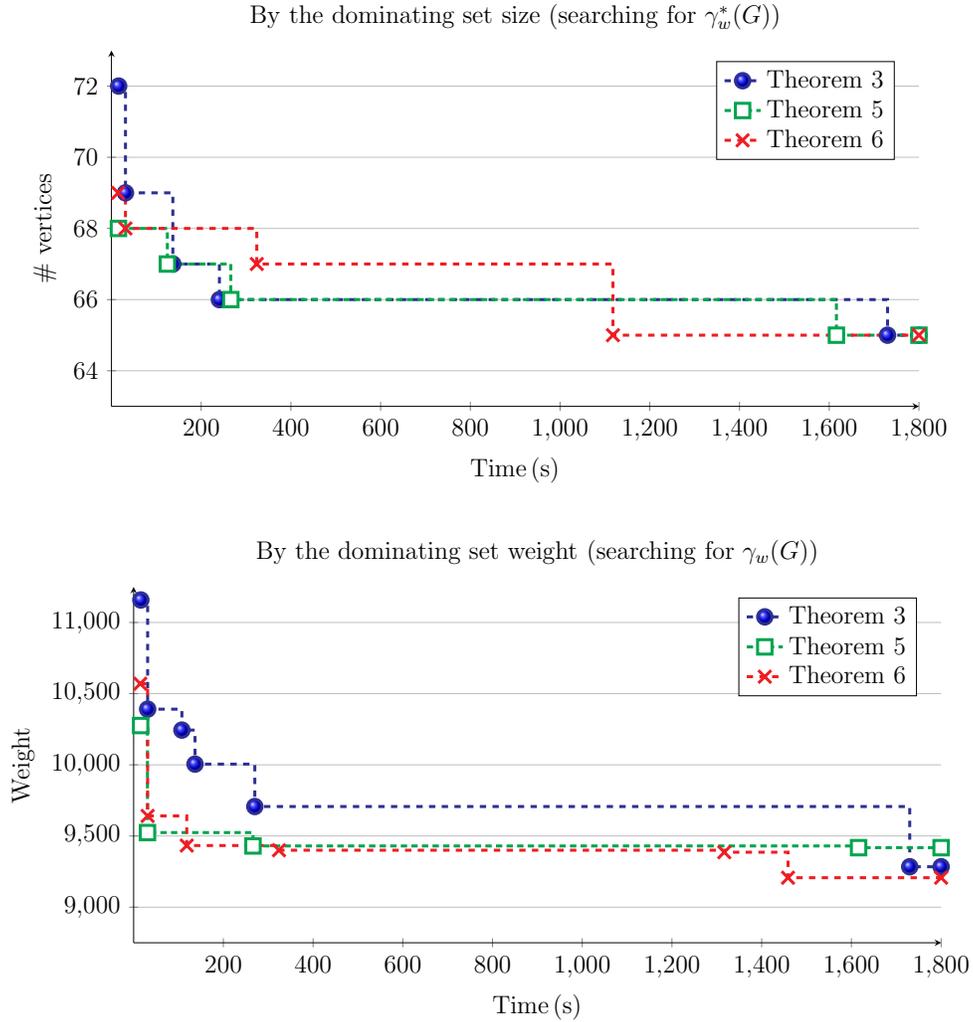

For the Erd\H{o}s--R\'enyi graph on $n=40,000$ veritces, the randomized algorithms arising from Theorems~\ref{uniform}, \ref{non-uniform}, and \ref{thm-inverse} respectively used $129$, $132$, and $132$ iterations.
For this graph, the ILP generic solver Gurobi ran out of memory and was not able to produce any solution. 
The performance profiles of the three randomized algorithms on this large size graph are presented in Figures~\ref{fig-ER-large-40-S} and Table~\ref{table-ER-large-40}.


\begin{table}[h!]
    \centering \footnotesize 
    \begin{tabular}{|c||ccc|ccc|ccc|}
    \hline
    \multirow{2}{*}{Problem} & \multicolumn{3}{c|}{Theorem~\ref{uniform}}                                                                                 & \multicolumn{3}{c|}{Theorem~\ref{non-uniform}}                                                                                 & \multicolumn{3}{c|}{Theorem~\ref{thm-inverse}}                                                                                 \\ \cline{2-10} 
                           & \multicolumn{1}{c|}{Size} & \multicolumn{1}{c|}{Wt} & \begin{tabular}[c]{@{}c@{}}Time\\ (sec)\end{tabular} & \multicolumn{1}{c|}{Size} & \multicolumn{1}{c|}{Wt} & \begin{tabular}[c]{@{}c@{}}Time\\ (sec)\end{tabular} & \multicolumn{1}{c|}{Size} & \multicolumn{1}{c|}{Wt} & \begin{tabular}[c]{@{}c@{}}Time\\ (sec)\end{tabular} \\ \hline \hline
    $\gamma^*_w(G)$               & \multicolumn{1}{c|}{65}   & \multicolumn{1}{c|}{9284}   & 1730                                                 & \multicolumn{1}{c|}{65}   & \multicolumn{1}{c|}{9418}   & 1616                                                 & \multicolumn{1}{c|}{65}   & \multicolumn{1}{c|}{9654}   & 1118                                                 \\ \hline
    $\gamma_w(G)$               & \multicolumn{1}{c|}{65}   & \multicolumn{1}{c|}{9284}   & 1730                                                 & \multicolumn{1}{c|}{65}   & \multicolumn{1}{c|}{9418}   & 1616                                                 & \multicolumn{1}{c|}{67}   & \multicolumn{1}{c|}{9207}   & 1459                                                 \\ \hline
    \end{tabular}
    \caption{\label{table-ER-large-40}Best found results for the Erd\H{o}s--R\'enyi graph on $40,000$ veritces.}
\end{table}

The randomized algorithms find the same quality solutions by the dominating set size when trying to find $\gamma^*_w(G)$, with the randomized algorithm of Theorem~\ref{thm-inverse} finding a set of this size quicker than the other two. However, the randomized algorithm of Theorem~\ref{uniform} finds the best (out of three) dominating set by weight of the same size. 
When searching for $\gamma_w(G)$, the best (out of three) solutions is found by the randomized algorithm of Theorem~\ref{thm-inverse}.
However, this dominating set is two vertices larger than the dominating sets found by the randomized algorithms of Theorem~\ref{uniform} and Theorem~\ref{non-uniform}.


\subsection{Sun graphs}
\label{sec-sun}
Given a minimum vertex degree $\delta\in \mathbb{N}$, $\delta \ge 3$, we define the random \emph{sun graph} model, denoted by $SG(\delta)$, following the description in \cite{ZP2011}. Assuming we have a complete graph $K_{\lfloor \delta\ln \delta \rfloor}$ with the vertex set $V(K_{\lfloor \delta\ln \delta \rfloor})$, we add an independent set of $\delta$ other vertices $O_{\delta}=\{v_1,v_2,...,v_{\delta}\}$ in such a way that each vertex $v_i$, $i=1,2,...,\delta$, is adjacent to exactly $\delta$ vertices of the set $V(K_{\lfloor \delta\ln \delta \rfloor})$, chosen uniformly at random. The resulting graph $G\in SG(\delta)$ has the vertex set $V(G)=V(K_{\lfloor \delta\ln \delta \rfloor})\cup O_{\delta}$ and $m={{\lfloor \delta\ln \delta \rfloor}\choose{2}} + \delta^2$ edges.
We call the resulting graph $G$ a \emph{sun graph}, which is clearly a split graph on $n=\lfloor \delta\ln \delta \rfloor + \delta$ vertices with $\delta(G)=\delta$.


Random (unweighted) sun graphs were generated for $\delta =50k$, $k=1,2,...,7$, ten graphs for each $k$, seventy graphs in total, {to form the set $SG(\delta=50k)$}.
Then, similarly to the case of Erd\H{o}s--R\'enyi graphs, for each $G\in SG(\delta=50k)$, $k=1,2,...,7$, weights were assigned to its vertices as integers in the range 
$\{101,102,..., 200\}$, uniformly at random.

\subsubsection{{Randomized heuristics and sun graphs}}
\label{sun-random}

Three randomized heuristics were run on each of these seventy sun graph instances. The heuristics correspond to Algorithm \ref{alg1}, using the probabilistic constructions of Theorems \ref{uniform}, \ref{non-uniform}, and \ref{thm-inverse}, respectively, 
The aggregated averages for ten sun graph instances for each $\delta=50k$, $k=1,2,...,7$, are presented in {Table \ref{table-SG}, together with corresponding CPU run-times}.

\begin{table}[h!]
    \centering \footnotesize 
    \begin{tabular}{|c||ccc||ccc||ccc|}
    \hline
                                                                                     & \multicolumn{3}{c|}{Theorem \ref{uniform}}                                                       & \multicolumn{3}{c|}{Theorem \ref{non-uniform}}                                                                                                       & \multicolumn{3}{c|}{Theorem \ref{thm-inverse}}                                                                                                       \\ \cline{2-10} 
    \multirow{-2}{*}{\begin{tabular}[c]{@{}c@{}}Min degree,\\ $\delta=50k$\end{tabular}} & \multicolumn{1}{c|}{Size} & \multicolumn{1}{c|}{Wt} & Time(s)                    & \multicolumn{1}{c|}{Size}                         & \multicolumn{1}{c|}{Wt}                         & Time(s)                    & \multicolumn{1}{c|}{Size}                         & \multicolumn{1}{c|}{Wt}                         & Time(s)                    \\ \hline \hline
                                                                                     & \multicolumn{1}{c|}{6.1}  & \multicolumn{1}{c|}{876}    &                            & \multicolumn{1}{c|}{\cellcolor[HTML]{FFCE93}5.8}  & \multicolumn{1}{c|}{783.9}                          &                            & \multicolumn{1}{c|}{6}                            & \multicolumn{1}{c|}{822.5}                          &                            \\ \cline{2-3} \cline{5-6} \cline{8-9}
    \multirow{-2}{*}{$k=1$}                                                           & \multicolumn{1}{c|}{6.3}  & \multicolumn{1}{c|}{864}    & \multirow{-2}{*}{0.11} & \multicolumn{1}{c|}{5.8}                          & \multicolumn{1}{c|}{\cellcolor[HTML]{FFCE93}783.9}  & \multirow{-2}{*}{0.12} & \multicolumn{1}{c|}{6.1}                          & \multicolumn{1}{c|}{812.8}                          & \multirow{-2}{*}{0.11} \\ \hline
                                                                                     & \multicolumn{1}{c|}{9}    & \multicolumn{1}{c|}{1279.7} &                            & \multicolumn{1}{c|}{8.7}                          & \multicolumn{1}{c|}{1221.6}                         &                            & \multicolumn{1}{c|}{\cellcolor[HTML]{FFCE93}8.5}  & \multicolumn{1}{c|}{1164.5}                         &                            \\ \cline{2-3} \cline{5-6} \cline{8-9}
    \multirow{-2}{*}{$k=2$}                                                           & \multicolumn{1}{c|}{9}    & \multicolumn{1}{c|}{1279.7} & \multirow{-2}{*}{0.43}  & \multicolumn{1}{c|}{9.1}                          & \multicolumn{1}{c|}{1210.2}                         & \multirow{-2}{*}{0.44}  & \multicolumn{1}{c|}{8.6}                          & \multicolumn{1}{c|}{\cellcolor[HTML]{FFCE93}1155.3} & \multirow{-2}{*}{0.43}  \\ \hline
                                                                                     & \multicolumn{1}{c|}{11}   & \multicolumn{1}{c|}{1669}   &                            & \multicolumn{1}{c|}{10.9}                         & \multicolumn{1}{c|}{1503.6}                         &                            & \multicolumn{1}{c|}{\cellcolor[HTML]{FFCE93}10.8} & \multicolumn{1}{c|}{1457}                           &                            \\ \cline{2-3} \cline{5-6} \cline{8-9}
    \multirow{-2}{*}{$k=3$}                                                         & \multicolumn{1}{c|}{11.3} & \multicolumn{1}{c|}{1645.8} & \multirow{-2}{*}{1.04}  & \multicolumn{1}{c|}{11.2}                         & \multicolumn{1}{c|}{1493.6}                         & \multirow{-2}{*}{1.04}  & \multicolumn{1}{c|}{10.9}                         & \multicolumn{1}{c|}{\cellcolor[HTML]{FFCE93}1456.2} & \multirow{-2}{*}{1.03}  \\ \hline
                                                                                     & \multicolumn{1}{c|}{12.5} & \multicolumn{1}{c|}{1818.3} &                            & \multicolumn{1}{c|}{12.3}                         & \multicolumn{1}{c|}{1671.5}                         &                            & \multicolumn{1}{c|}{\cellcolor[HTML]{FFCE93}12.2} & \multicolumn{1}{c|}{1734.6}                         &                            \\ \cline{2-3} \cline{5-6} \cline{8-9}
    \multirow{-2}{*}{$k=4$}                                                          & \multicolumn{1}{c|}{12.8} & \multicolumn{1}{c|}{1802}   & \multirow{-2}{*}{2.09}   & \multicolumn{1}{c|}{12.5}                         & \multicolumn{1}{c|}{\cellcolor[HTML]{FFCE93}1671.2} & \multirow{-2}{*}{2.14}   & \multicolumn{1}{c|}{12.5}                         & \multicolumn{1}{c|}{1681.3}                         & \multirow{-2}{*}{2.14}   \\ \hline
                                                                                     & \multicolumn{1}{c|}{14}   & \multicolumn{1}{c|}{2070.4} &                            & \multicolumn{1}{c|}{\cellcolor[HTML]{FFCE93}13.7} & \multicolumn{1}{c|}{1905.8}                         &                            & \multicolumn{1}{c|}{13.8}                         & \multicolumn{1}{c|}{1959}                           &                            \\ \cline{2-3} \cline{5-6} \cline{8-9}
    \multirow{-2}{*}{$k=5$}                                                           & \multicolumn{1}{c|}{14.3} & \multicolumn{1}{c|}{2048.3} & \multirow{-2}{*}{3.09}   & \multicolumn{1}{c|}{14}                           & \multicolumn{1}{c|}{\cellcolor[HTML]{FFCE93}1872.2} & \multirow{-2}{*}{3.11}   & \multicolumn{1}{c|}{14.2}                         & \multicolumn{1}{c|}{1927.4}                         & \multirow{-2}{*}{3.07}     \\ \hline
                                                                                     & \multicolumn{1}{c|}{15.1} & \multicolumn{1}{c|}{2211.6} &                            & \multicolumn{1}{c|}{15.1}                         & \multicolumn{1}{c|}{2079.3}                         &                            & \multicolumn{1}{c|}{\cellcolor[HTML]{FFCE93}14.9} & \multicolumn{1}{c|}{2143.5}                         &                            \\ \cline{2-3} \cline{5-6} \cline{8-9}
    \multirow{-2}{*}{$k=6$}                                                          & \multicolumn{1}{c|}{15.3} & \multicolumn{1}{c|}{2194.6} & \multirow{-2}{*}{4.95}   & \multicolumn{1}{c|}{15.3}                         & \multicolumn{1}{c|}{\cellcolor[HTML]{FFCE93}2061.3} & \multirow{-2}{*}{4.8}   & \multicolumn{1}{c|}{15.4}                         & \multicolumn{1}{c|}{2088.7}                         & \multirow{-2}{*}{4.76}   \\ \hline
                                                                                     & \multicolumn{1}{c|}{16}   & \multicolumn{1}{c|}{2375}   &                            & \multicolumn{1}{c|}{\cellcolor[HTML]{FFCE93}15.5} & \multicolumn{1}{c|}{2147.5}                         &                            & \multicolumn{1}{c|}{16}                           & \multicolumn{1}{c|}{2245.2}                         &                            \\ \cline{2-3} \cline{5-6} \cline{8-9}
    \multirow{-2}{*}{$k=7$}                                                          & \multicolumn{1}{c|}{16.1} & \multicolumn{1}{c|}{2371.7} & \multirow{-2}{*}{6.56}   & \multicolumn{1}{c|}{15.6}                         & \multicolumn{1}{c|}{\cellcolor[HTML]{FFCE93}2145.5} & \multirow{-2}{*}{6.6}   & \multicolumn{1}{c|}{16.2}                         & \multicolumn{1}{c|}{2213.1}                         & \multirow{-2}{*}{6.62}    \\ \hline
    \end{tabular}
    \caption{\label{table-SG}Aggregated results of running the randomized heuristics on the sun graphs.}
\end{table}

Similarly to the case of Erd\H{o}s--R\'enyi graphs, each of the randomized algorithms was run {twenty} times on each sun graph instance, and the best found solution (out of {twenty}) by the dominating set size and also by the dominating set weight were recorded in {Table \ref{table-SG}} (in the upper and lower subrows for each $k=1,...,7$, respectively). 
{Table \ref{table-SG}} shows that the randomized algorithm corresponding to the probabilistic construction of Theorem \ref{uniform} performs worse on sun graphs than those corresponding to Theorems \ref{non-uniform} and \ref{thm-inverse}. 
The randomized algorithms corresponding to Theorems \ref{uniform}, \ref{non-uniform}, and \ref{thm-inverse}  find better dominating sets by size for {$52.9\%, 71.4\%$, and $70\%$} of all instances, respectively, and better dominating sets by weight for {$7.1\%, 42.9\%$, and $50\%$} of all instances, respectively.

The upper bounds of Theorems \ref{Caro} and \ref{uniform} were satisfied for all the problem instances.
{Table \ref{table-SG}} also shows that the heuristic methods derived from Theorems \ref{non-uniform} and \ref{thm-inverse} provide better results for the two-objective optimization problem than that of Theorem \ref{uniform} in all the cases, although both methods of Theorems \ref{non-uniform} and \ref{thm-inverse} are designed for optimization by weight only.

\subsubsection{{Using a generic solver on sun graphs}}
\label{sun-generic}

For each weighted sun graph $G\in SG(\delta=50k)$, $k=1,2,...,7$, we attempted to find exact deterministic solutions to the problems of finding $\gamma^*_w(G)$ and $\gamma_w(G)$ by using the ILP formulation (\ref{ILP1}) and 
{the ILP generic solver Gurobi \cite{Gurobi}}.
This was successful in a reasonable amount of {CPU time of at most $30$ minutes ($1800$ sec)} only for $k=1,2$ ($245$ and $560$ vertices, respectively) {when computing $\gamma^*_w(G)$, and for $k=1,2,3$ ($245, 560,$ and $901$ vertices, respectively) when computing $\gamma_w$}. 
In these deterministic computational experiments, almost half of the graph instances (nine out of twenty) were found to have a dominating set of a non-minimum size with weight lower than $\gamma^*_w(G)$, i.e. $\gamma_w(G)<\gamma^*_w(G)$ (see {Table \ref{table-sun-Gurobi}}).

Similarly to the Erd\H{o}s--R\'enyi graphs, for the larger (medium) size sun graph instances, i.e. for $k\ge 3$, when computing $\gamma^*_w(G)$, and for $k\ge 4$, when computing $\gamma_w(G)$, Gurobi was run for $30$ minutes ($1800$ sec CPU time) as a benchmark heuristic solver. The aggregated results are presented in Table~\ref{table-sun-Gurobi} (the deterministic solution results are shaded).

\begin{table}[H]
    \centering \footnotesize 
    \begin{tabular}{|c||ccccc||ccccc|}
    \hline
    \multirow{3}{*}{\begin{tabular}[c]{@{}c@{}}Min\\ degree,\\ $\delta=50k$\end{tabular}}  & \multicolumn{5}{c|}{ILP for $\gamma^*_w(G)$}                                                                                                                                                                                                                                                        & \multicolumn{5}{c|}{ILP for $\gamma_w(G)$}                                                                                                                                                                                                                                                       \\ \cline{2-11} 
                                                                                       & \multicolumn{1}{c|}{\begin{tabular}[c]{@{}c@{}}IHS\\ size\end{tabular}} & \multicolumn{1}{c|}{\begin{tabular}[c]{@{}c@{}}IHS CPU\\ time\,(s)\end{tabular}} & \multicolumn{1}{c|}{\begin{tabular}[c]{@{}c@{}}BPS\\ size\end{tabular}} &  \multicolumn{1}{c|}{\begin{tabular}[c]{@{}c@{}}BPS\\ wt\end{tabular}} &  \multicolumn{1}{c|}{\begin{tabular}[c]{@{}c@{}}BPS CPU\\ time\,(s)\end{tabular}}   &  \multicolumn{1}{c|}{\begin{tabular}[c]{@{}c@{}}IHS\\ size\end{tabular}}  & \multicolumn{1}{c|}{\begin{tabular}[c]{@{}c@{}}IHS CPU\\ time\,(s)\end{tabular}} & \multicolumn{1}{c|}{\begin{tabular}[c]{@{}c@{}}BPS\\ wt\end{tabular}} & \multicolumn{1}{c|}{\begin{tabular}[c]{@{}c@{}}BPS\\ size\end{tabular}} & \multicolumn{1}{c|}{\begin{tabular}[c]{@{}c@{}}BPS CPU\\ time\,(s)\end{tabular}}  \\ \hline \hline
    $k=1$                                                                                & \multicolumn{1}{c|}{50}                                                       & \multicolumn{1}{c|}{0.026}                                                  & \multicolumn{1}{c|}{\cellcolor{orange!25}4.5}                                                    & \multicolumn{1}{c|}{\cellcolor{orange!25}583.1}  & \cellcolor{orange!25}2.78 & \multicolumn{1}{c|}{7623.5}                                                   & \multicolumn{1}{c|}{0.031}                                                  & \multicolumn{1}{c|}{\cellcolor{orange!25}564.9}                                                    & \multicolumn{1}{c|}{\cellcolor{orange!25}4.9}  & \cellcolor{orange!25}0.9  \\ \hline
    $k=2$                                                                                & \multicolumn{1}{c|}{100}                                                      & \multicolumn{1}{c|}{0.18}                                                  & \multicolumn{1}{c|}{\cellcolor{orange!25}6}                                                      & \multicolumn{1}{c|}{\cellcolor{orange!25}758.8}  & \cellcolor{orange!25}183.68 & \multicolumn{1}{c|}{15175.6}                                                  & \multicolumn{1}{c|}{0.19}                                                  & \multicolumn{1}{c|}{\cellcolor{orange!25}747.1}                                                    & \multicolumn{1}{c|}{\cellcolor{orange!25}6.5}  & \cellcolor{orange!25}31.44  \\ \hline
    $k=3$                                                                                & \multicolumn{1}{c|}{150}                                                      & \multicolumn{1}{c|}{0.47}                                                  & \multicolumn{1}{c|}{7.8}                                                    & \multicolumn{1}{c|}{947.7}  & 1800   & \multicolumn{1}{c|}{22540.9}                                                  & \multicolumn{1}{c|}{0.46}                                                  & \multicolumn{1}{c|}{\cellcolor{orange!25}859.8}                                                    & \multicolumn{1}{c|}{\cellcolor{orange!25}8}    & \cellcolor{orange!25}1381.21 \\ \hline
    $k=4$                                                                                & \multicolumn{1}{c|}{200}                                                      & \multicolumn{1}{c|}{0.92}                                                  & \multicolumn{1}{c|}{8.8}                                                    & \multicolumn{1}{c|}{1120}   & 1800   & \multicolumn{1}{c|}{30105.4}                                                  & \multicolumn{1}{c|}{0.91}                                                  & \multicolumn{1}{c|}{1005.9}                                                   & \multicolumn{1}{c|}{9.2}  & 1800    \\ \hline
    $k=5$                                                                                & \multicolumn{1}{c|}{250}                                                      & \multicolumn{1}{c|}{1.57}                                                  & \multicolumn{1}{c|}{9.9}                                                    & \multicolumn{1}{c|}{1169.1} & 1800   & \multicolumn{1}{c|}{37734}                                                    & \multicolumn{1}{c|}{1.56}                                                  & \multicolumn{1}{c|}{1122.5}                                                   & \multicolumn{1}{c|}{10.2} & 1800    \\ \hline
    $k=6$                                                                                & \multicolumn{1}{c|}{300}                                                      & \multicolumn{1}{c|}{2.51}                                                  & \multicolumn{1}{c|}{10.8}                                                   & \multicolumn{1}{c|}{1341.3} & 1800   & \multicolumn{1}{c|}{45061.2}                                                  & \multicolumn{1}{c|}{2.48}                                                  & \multicolumn{1}{c|}{1218.3}                                                   & \multicolumn{1}{c|}{11.2} & 1800    \\ \hline
    $k=7$                                                                                & \multicolumn{1}{c|}{350}                                                      & \multicolumn{1}{c|}{3.76}                                                   & \multicolumn{1}{c|}{11.4}                                                   & \multicolumn{1}{c|}{1489.2} & 1800   & \multicolumn{1}{c|}{52535}                                                    & \multicolumn{1}{c|}{3.71}                                                  & \multicolumn{1}{c|}{1283.3}                                                   & \multicolumn{1}{c|}{11.9} & 1800    \\ \hline
    \end{tabular}
    \caption{\label{table-sun-Gurobi}Aggregated results of running the ILP generic solver on the sun graphs.}
\end{table}

It can be seen from Tables~\ref{table-SG} and \ref{table-sun-Gurobi} that, for the small and medium size sun graphs, the new randomized heuristics provide solutions of about one order of magnitude better than the IHS of Gurobi, but Gurobi finds its (apparently trivial) IHS faster.
It can also be seen that, on average, the best solutions found by the randomized algorithms for {$k=1,2,3$} are about {$39-69\%$} worse by weight and, {for $k=1,2$,} about {$29-42\%$}  worse by size than the optimal {(deterministic)} solutions  (see shaded cells in BPS columns in Table~\ref{table-sun-Gurobi}). 
The computing time of the deterministic solution methods (ILP) becomes prohibitively high for graphs of {$1259$ and more} vertices in these experiments. 
{In comparison to the heuristic results for medium size sun graphs ($k\ge 4$) obtained by Gurobi in 1800 sec, the very quick (less than ten seconds) heuristic solutions by the new randomized algorithms are about $66-69\%$ worse by weight and about $36-39\%$ worse by size. Again, this gap can be significantly reduced and even better results can be obtained by running the randomized algorithms for more iterations for the same amount of CPU time ($1800$ sec). In Section~\ref{sun-large}, we show that the randomized algorithms provide clearly better results for large size sun graphs when run for the same amount of time as Gurobi.}

{Similarly to the Erd\H{o}s--R\'enyi graphs}, the average CPU run-times to compute exact values of $\gamma^*_w(G)$ and $\gamma_w(G)$ (for $k=1,2,3$) using the ILP formulations and {generic solver 
Gurobi (see BPS columns in Table~\ref{table-sun-Gurobi})} clearly show the considerable growth of computational time requirements with respect to the graph order, as well as differences in computational complexity for finding $\gamma^*_w(G)$ and $\gamma_w(G)$ in sun graphs.
It can be seen from {Tables \ref{table-ER-Gurobi} and \ref{table-sun-Gurobi}} that running {Gurobi} on the ILP formulations of the problems is much more efficient for the sun graphs than for the Erd\H{o}s--R\'enyi graphs. This may be explained by the structure of the sun graphs, which are split graphs. However, both types of random graphs exhibit exponential growth in the ILP solution time with respect to the graph order. 

\subsubsection{{Large-scale sun graphs}}
\label{sun-large}

Finally, we have run the three randomized algorithms and the ILP generic solver Gurobi on two large size sun graphs, generated as described at the beginning of  Section~\ref{sec-sun}. 
One graph has the minimum vertex degree $\delta=1250$ ($k=25$) and $n=10163$ vertices, the other has $\delta=2500$ ($k=50$) and $n=22060$ vertices. 
Each of the four solvers was given $30$ minutes (1800 sec) of CPU time to find a solution to the problems corresponding to searching for $\gamma_w^*(G)$ and $\gamma_w(G)$.

For the sun graph on $10163$ vertices ($\delta=1250$), the randomized algorithms arising from Theorems~\ref{uniform}, \ref{non-uniform}, and \ref{thm-inverse} respectively used  $125$, $101$, and $87$ iterations.
For this graph, Gurobi was able to find only the initial trivial solution of $1250$ vertices (in $128$ sec) and weight $189407$ (in $116$ sec). No other solution was found by the ILP generic solver in $1800$ sec.
The performance profiles of the randomized algorithms on this large sun graph are presented in Figure~\ref{fig-sun-large-1} and Table~\ref{table-sun-large-1} below.

\begin{figure}[H]
    \centering
    \begin{tikzpicture}[scale=0.8]
        \begin{axis}[
            axis lines=left,
            title={By the dominating set size\  (searching for $\gamma^*_w(G)$)},
            xlabel={Time\,(s)},
            ylabel={\# vertices},
            xmin=0, xmax=1800,
            ymin=22, ymax=32,
            xtick={200,400,600,800,1000,1200,1400,1600,1800},
            ytick={23,25,27,29,31},
            ymajorgrids=true,
            legend pos=north east,
        ]
        
        \addplot+[
            const plot,
            color=Blue,
            mark=ball,
            mark options={fill=Blue,style=solid},
            mark size=1.2mm,
            style=dashed,
            line width=0.5mm
            ] coordinates {
                (12.1846,30)(37.7752,28)(182.699,26)(530.244,25)(1245.96,24)(1800,24)
            };
            \addlegendentry{Theorem 3}
        
        \addplot+[
            const plot,
            color=Green,
            mark=square*,
            mark options={fill=White,style=solid},
            mark size=1.2mm,
            style=densely dashed,
            line width=0.5mm
            ] coordinates {
                (14.5633,29)(31.0361,25)(1800,25)
            };
            \addlegendentry{Theorem 5}

        \addplot+[
            const plot,
            color=Red,
            mark=x,
            mark options={fill=Red,style=solid},
            mark size=1.5mm,
            style=dashed,
            line width=0.5mm
            ] coordinates {
                (16.4106,28)(58.9779,27)(251.044,25)(1800,25)
            };
            \addlegendentry{Theorem 6}
        \end{axis}
    \end{tikzpicture}

\vspace{5mm}

    \begin{tikzpicture}[scale=0.8]
        \begin{axis}[
            axis lines=left,
            title={By the dominating set weight\  (searching for $\gamma_w(G)$)},
            xlabel={Time\,(s)},
            ylabel={Weight},
            xmin=0, xmax=1800,
            ymin=3400, ymax=4800,
            xtick={200,400,600,800,1000,1200,1400,1600,1800},
            ytick={3500,3750,4000,4250,4500,4750},
            ymajorgrids=true,
            legend pos=north east,
        ]
        
        \addplot+[
            const plot,
            color=Blue,
            mark=ball,
            mark options={fill=Blue,style=solid},
            mark size=1.4mm,
            style=dashed,
            line width=0.5mm
            ] coordinates {
                (12.1846,4800)(25.2619,4669)(37.7752,4206)(127.116,4135)(369.995,4050)(530.244,3805)(1245.96,3557)(1800,3557)
            };
            \addlegendentry{Theorem 3}
        
        \addplot+[
            const plot,
            color=Green,
            mark=square*,
            mark options={fill=White,style=solid},
            mark size=1mm,
            style=densely dashed,
            line width=0.5mm
            ] coordinates {
                (14.5633,4417)(31.0361,3736)(284.425,3613)(755.876,3558)(1800,3558)
            };
            \addlegendentry{Theorem 5}

        \addplot+[
            const plot,
            color=Red,
            mark=x,
            mark options={fill=Red,style=solid},
            mark size=1.5mm,
            style=dashed,
            line width=0.5mm
            ] coordinates {
                (16.4106,3870)(143.414,3859)(165.446,3755)(251.044,3461)(1800,3461)
            };
            \addlegendentry{Theorem 6}
        \end{axis}
    \end{tikzpicture}
     \caption{\label{fig-sun-large-1}{Running the randomized heuristics on the sun graph of $10163$ vertices ($\delta=1250$).}}
\end{figure}
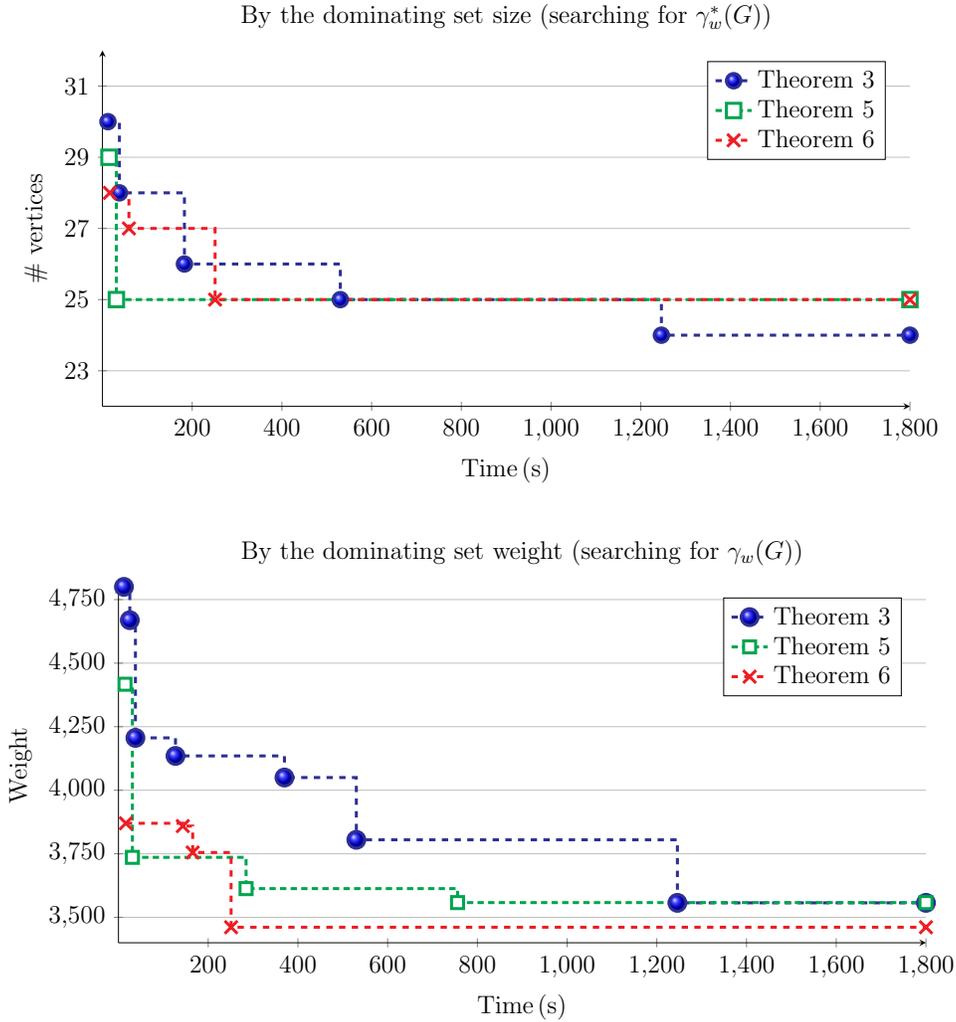


\begin{table}[H]
    \centering \footnotesize
    {
    \begin{tabular}{|c||ccc||ccc||ccc|}
    \hline 
    \multirow{2}{*}{Problem} & \multicolumn{3}{c|}{Theorem~\ref{uniform}}                                                                                 & \multicolumn{3}{c|}{Theorem~\ref{non-uniform}}                                                                                 & \multicolumn{3}{c|}{Theorem~\ref{thm-inverse}}                                                                                 \\ \cline{2-10} 
                           & \multicolumn{1}{c|}{Size} & \multicolumn{1}{c|}{Wt} & \begin{tabular}[c]{@{}c@{}}Time\\ (sec)\end{tabular} & \multicolumn{1}{c|}{Size} & \multicolumn{1}{c|}{Wt} & \begin{tabular}[c]{@{}c@{}}Time\\ (sec)\end{tabular} & \multicolumn{1}{c|}{Size} & \multicolumn{1}{c|}{Wt} & \begin{tabular}[c]{@{}c@{}}Time\\ (sec)\end{tabular} \\ \hline
    $\gamma^*_w(G)$                   & \multicolumn{1}{c|}{24}   & \multicolumn{1}{c|}{3557}   & 1246                                                 & \multicolumn{1}{c|}{25}   & \multicolumn{1}{c|}{3726}   & 1259                                                 & \multicolumn{1}{c|}{25}   & \multicolumn{1}{c|}{3461}   & 251                                                  \\ \hline
    $\gamma_w(G)$                 & \multicolumn{1}{c|}{24}   & \multicolumn{1}{c|}{3557}   & 1246                                                 & \multicolumn{1}{c|}{27}   & \multicolumn{1}{c|}{3558}   & 756                                                  & \multicolumn{1}{c|}{25}   & \multicolumn{1}{c|}{3461}   & 251                                                  \\ \hline
    \end{tabular}
    \caption{\label{table-sun-large-1}Best found results for the sun graph on $10163$ veritces.}
    }
\end{table}

When trying to find $\gamma^*_w(G)$, the randomized algorithm arising from Theorem~\ref{uniform} provides the best solution by the dominating set size. Again, this can be explained by the same optimal probability used in Theorems~\ref{uniform} and \ref{Caro}. 
When searching for $\gamma_w(G)$, the randomized algorithm arising from Theorem~\ref{thm-inverse} provides the best solution, which, similarly to the large Erd\H{o}s--R\'enyi graphs, can be explained by its sensitivity to the vertex weights in the graph.

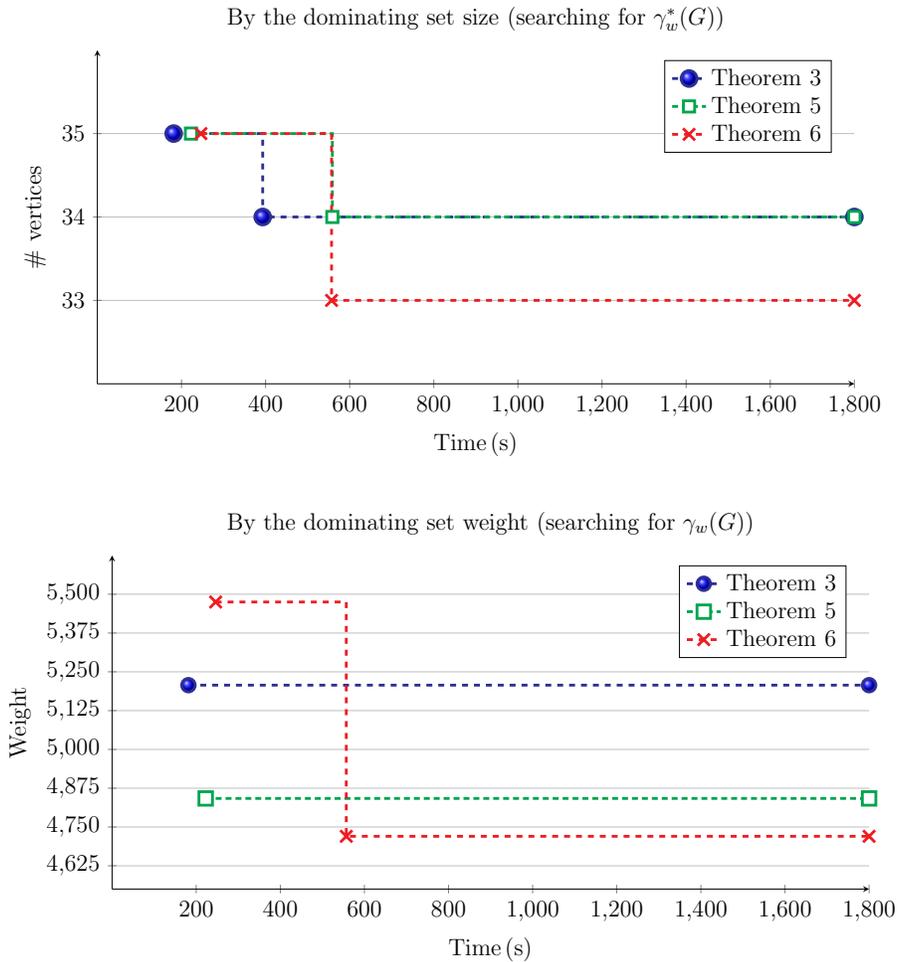
\begin{figure}[H]
    \centering
    \begin{tikzpicture}[scale=0.75]
        \begin{axis}[
            axis lines=left,
            title={By the dominating set size\  (searching for $\gamma^*_w(G)$)},
            xlabel={Time\,(s)},
            ylabel={\# vertices},
            xmin=0, xmax=1800,
            ymin=32, ymax=36,
            xtick={200,400,600,800,1000,1200,1400,1600,1800},
            ytick={33,34,35},
            ymajorgrids=true,
            legend pos=north east,
        ]
        
        \addplot+[
            const plot,
            color=Blue,
            mark=ball,
            mark options={fill=Blue,style=solid},
            mark size=1.4mm,
            style=dashed,
            line width=0.5mm
            ] coordinates {
                (181.523,35)(393.174,34)(1800,34)
            };
            \addlegendentry{Theorem 3}
        
        \addplot+[
            const plot,
            color=Green,
            mark=square*,
            mark options={fill=White,style=solid},
            mark size=1mm,
            style=densely dashed,
            line width=0.5mm
            ] coordinates {
                (222.627,35)(558.807,34)(1800,34)
            };
            \addlegendentry{Theorem 5}

        \addplot+[
            const plot,
            color=Red,
            mark=x,
            mark options={fill=Red,style=solid},
            mark size=1.5mm,
            style=dashed,
            line width=0.5mm
            ] coordinates {
                (246.274,35)(556.917,33)(1800,33)
            };
            \addlegendentry{Theorem 6}
        \end{axis}
    \end{tikzpicture}

\vspace{5mm}

    \begin{tikzpicture}[scale=0.75]
        \begin{axis}[
            axis lines=left,
            title={By the dominating set weight\  (searching for $\gamma_w(G)$)},
            xlabel={Time\,(s)},
            ylabel={Weight},
            xmin=0, xmax=1800,
            ymin=4550, ymax=5625,
            xtick={200,400,600,800,1000,1200,1400,1600,1800},
            ytick={4625,4750,4875,5000,5125,5250,5375,5500},
            ymajorgrids=true,
            legend pos=north east,
        ]
        
        \addplot+[
            const plot,
            color=Blue,
            mark=ball,
            mark options={fill=Blue,style=solid},
            mark size=1.2mm,
            style=dashed,
            line width=0.5mm
            ] coordinates {
                (181.523,5207)(1800,5207)
            };
            \addlegendentry{Theorem 3}
        
        \addplot+[
            const plot,
            color=Green,
            mark=square*,
            mark options={fill=White,style=solid},
            mark size=1.2mm,
            style=densely dashed,
            line width=0.5mm
            ] coordinates {
                (222.627,4842)(1800,4842)
            };
            \addlegendentry{Theorem 5}

        \addplot+[
            const plot,
            color=Red,
            mark=x,
            mark options={fill=Red,style=solid},
            mark size=1.5mm,
            style=dashed,
            line width=0.5mm
            ] coordinates {
                (246.274,5475)(556.917,4720)(1800,4720)
            };
            \addlegendentry{Theorem 6}
        \end{axis}
    \end{tikzpicture}
    \caption{\label{fig-sun-large-2}{Running the randomized heuristics on the sun graph of $22060$ vertices ($\delta=2500$).}}
\end{figure}

For the sun graph on $22060$ vertices ($\delta=2500$), the randomized algorithms arising from Theorems~\ref{uniform}, \ref{non-uniform}, and \ref{thm-inverse} respectively used  $11$, $9$, and $8$ iterations.
For this graph, Gurobi ran out of memory and was not able to produce any solution. 
The performance profiles of the randomized algorithms on this large sun graph are presented in Figure~\ref{fig-sun-large-2} and Table~\ref{table-sun-large-2}.

\begin{table}[H]
    \centering \footnotesize
    {
    \begin{tabular}{|c||ccc||ccc||ccc|}
    \hline
    \multirow{2}{*}{Problem} & \multicolumn{3}{c|}{Theorem~\ref{uniform}}                                                                                 & \multicolumn{3}{c|}{Theorem~\ref{non-uniform}}                                                                                 & \multicolumn{3}{c|}{Theorem~\ref{thm-inverse}}                                                                                 \\ \cline{2-10} 
                           & \multicolumn{1}{c|}{Size} & \multicolumn{1}{c|}{Wt} & \begin{tabular}[c]{@{}c@{}}Time\\ (sec)\end{tabular} & \multicolumn{1}{c|}{Size} & \multicolumn{1}{c|}{Wt} & \begin{tabular}[c]{@{}c@{}}Time\\ (sec)\end{tabular} & \multicolumn{1}{c|}{Size} & \multicolumn{1}{c|}{Wt} & \begin{tabular}[c]{@{}c@{}}Time\\ (sec)\end{tabular} \\ \hline
    $\gamma^*_w(G)$                   & \multicolumn{1}{c|}{34}   & \multicolumn{1}{c|}{5253}   & 393                                                  & \multicolumn{1}{c|}{34}   & \multicolumn{1}{c|}{5054}   & 559                                                  & \multicolumn{1}{c|}{33}   & \multicolumn{1}{c|}{4720}   & 557                                                  \\ \hline
    $\gamma_w(G)$                 & \multicolumn{1}{c|}{35}   & \multicolumn{1}{c|}{5207}   & 182                                                  & \multicolumn{1}{c|}{35}   & \multicolumn{1}{c|}{4842}   & 223                                                  & \multicolumn{1}{c|}{33}   & \multicolumn{1}{c|}{4720}   & 557                                                  \\ \hline
    \end{tabular}
    \caption{\label{table-sun-large-2}Best found results for the sun graph on $22060$ veritces.}
     }
\end{table}

{
For this large sun graph, the best solutions are found by the randomized algorithm arising from Theorem~\ref{thm-inverse}. While this seems to be normal for the search by weight, the result by size illustrates the randomized nature of the three solvers.
}


\section{Conclusions}
\label{sec:conclusion}

We have considered the problem of finding a minimum-size dominating set of the smallest weight and the problem of finding a minimum weight dominating set in a vertex-weighted simple graph. 
The two-objective problem is different from the classic single-objective problems of finding a smallest cardinality dominating set or a smallest weight dominating set in a graph,
{but could be seen as a particular case of one of them}. 
We have shown three generalizations of the upper bounds for the domination number using the probabilistic method. 
The probabilistic constructions used to prove the new bounds and theorems provide randomized heuristics to find some approximate solutions for the weighted domination problems. 
We have shown a connection between the problems in graphs via ILP formulations of the problems and provided a reduction from the two-objective problem to the minimum weight dominating set problem.

{The simple and efficient randomized heuristics presented in this paper can be used, for example, to quickly find a better initial solution in the local search and other heuristics presented in \cite{AV2018,BB2016,WCCY2018}, or to help a generic ILP solver like Gurobi \cite{Gurobi} to start an exhaustive search with a better and quickly found initial heuristic solution. The proposed randomized heuristics are very efficient in usage of computer memory and CPU time.}
Since the probability used in the proof of Theorem \ref{uniform} is the same as in the proof of Theorem \ref{Caro}, the corresponding randomized algorithm from the framework of Algorithm \ref{alg1} provides a certain optimality in finding a good dominating set not only by weight, but also by size.
On the other hand, the computational experiments show that the randomized algorithms arising from Theorems \ref{non-uniform} and \ref{thm-inverse} are more sensitive to the vertex weight parameters and usually provide better results than the randomized algorithm arising from Theorems \ref{uniform}, in particular, in the case of sun graphs of Section \ref{sec-sun}. 

The upper bound and probabilistic construction of Theorem \ref{uniform} are reminiscent and similar to those used in Theorem \ref{Caro} for unweighted graphs, and the assumptions of Theorem \ref{uniform} require the graph weights not to vary too much among the vertices. On the other hand, ILP reduction (\ref{ILP2}) from Section \ref{sec:deterministic-heuristic} guarantees that the reduced problem has all the vertex weights close to one and not varying too much. In this case, the bounds and probabilities of Theorems \ref{non-uniform} and \ref{thm-inverse} turn out to be close to those of Theorem \ref{uniform}, which makes all three of the theorems useful in the context of reduction.

As a direction for future work, it would be interesting to consider a relaxation of the problems to find dominating sets whose size and/or weight are within a certain limit from $\gamma(G)$ and/or $\gamma_{w}(G)$ 
in a graph, respectively, and a connection between $\gamma(G)$ and vertex weight parameters. Also, it would be interesting to devise other heuristic and deterministic solution methods for the problems. In particular, it would be useful to develop some direct deterministic solution algorithms, like the state-of-the-art deterministic algorithms for the minimum-size dominating sets  \cite{Wendy2, Wendy1}. 
Eventually, the upper bounds and outcomes of the randomized algorithms presented in this paper can be used to reduce the search space in smart exhaustive search algorithms, like backtracking and branch-and-bound, for the minimum-weight dominating sets.\\





\end{document}